\documentclass[11pt,reqno]{amsart}
\usepackage[margin=1in]{geometry}
\usepackage{microtype} 
\usepackage{amsmath,amssymb,amsthm,enumerate,booktabs,needspace}
\usepackage[colorlinks=true,linkcolor=blue,citecolor=blue,urlcolor=blue]{hyperref}

\usepackage{amsmath,tikz}
\usetikzlibrary{arrows.meta,calc}

\definecolor{emInk}{HTML}{26333D}
\definecolor{emP}{HTML}{006B78}
\definecolor{emQ}{HTML}{9B3146}
\definecolor{emPale}{HTML}{EFF3F5}
\definecolor{emGrey}{HTML}{7E8B93}
\definecolor{emCut}{HTML}{E4EDF0}
\tikzset{
  emEdge/.style={draw=emInk,line width=.85pt,line cap=round,line join=round},
  emP/.style={emEdge,draw=emP,line width=1.3pt},
  emQ/.style={emEdge,draw=emQ,line width=1.3pt},
  emChain/.style={emP,line width=1.65pt},
  emTunnel/.style={emEdge,line width=1pt},
  emVertex/.style={circle,draw=emInk,fill=white,line width=.7pt,
    inner sep=0pt,minimum size=3.5pt},
  emTag/.style={font=\small,inner sep=1pt,text=emInk},
  emGhost/.style={draw=emGrey,densely dashed,line width=.65pt},
  emCut/.style={draw=emInk,fill=emCut,line width=.75pt},
  emArrow/.style={-{Stealth[length=4pt,width=3pt]},draw=emGrey,line width=.65pt}
}

\newcommand{\EvenPlatStubs}[8]{%
  \draw[#1] (-1.5,.65)--(-1.5,1.15)
    ..controls(-1.5,1.60)and(-1.25,1.85)..(-1,1.85);
  \draw[#2] (-1,1.85)..controls(-.75,1.85)and(-.5,1.60)..(-.5,1.15)--(-.5,.65);
  \draw[#3] (.5,.65)--(.5,1.15)
    ..controls(.5,1.60)and(.75,1.85)..(1,1.85);
  \draw[#4] (1,1.85)..controls(1.25,1.85)and(1.5,1.60)..(1.5,1.15)--(1.5,.65);
  \draw[#5] (-1.5,-.65)--(-1.5,-1.15)
    ..controls(-1.5,-1.60)and(-1.25,-1.85)..(-1,-1.85);
  \draw[#6] (-1,-1.85)..controls(-.75,-1.85)and(-.5,-1.60)..(-.5,-1.15)--(-.5,-.65);
  \draw[#7] (.5,-.65)--(.5,-1.15)
    ..controls(.5,-1.60)and(.75,-1.85)..(1,-1.85);
  \draw[#8] (1,-1.85)..controls(1.25,-1.85)and(1.5,-1.60)..(1.5,-1.15)--(1.5,-.65);
}
\newcommand{\EvenPlatBox}{%
  \path[draw=emGrey,fill=emPale,line width=.65pt,rounded corners=1.5pt]
    (-1.72,-.65) rectangle (1.72,.65);
  \node[text=emInk,font=\large] at (0,0) {$B_r$};
}
\newcommand{\EvenPlatEdgeNames}[4]{%
  \foreach \x/\n in {-1.5/#1,-.5/#2,.5/#3,1.5/#4}
    \node[font=\scriptsize,text=emInk,inner sep=0pt] at (\x,.40) {$\n$};
}
\newcommand{\EvenUpperTunnel}[1]{%
  \draw[emTunnel] (-1,1.85)..controls(-.65,2.66)and(.65,2.66)..(1,1.85);
  \node[emTag,above] at (0,2.46) {$\tau_+ : #1$};
  \node[emVertex] at (-1,1.85) {};\node[emVertex] at (1,1.85) {};
}
\newcommand{\EvenLowerTunnel}[1]{%
  \draw[emTunnel] (-1,-1.85)..controls(-.65,-2.66)and(.65,-2.66)..(1,-1.85);
  \node[emTag,below] at (0,-2.46) {$\tau_- : #1$};
  \node[emVertex] at (-1,-1.85) {};\node[emVertex] at (1,-1.85) {};
}
\newcommand{\EvenTopCapLabels}[2]{%
  \node[emTag,above] at (-1,1.97) {$#1$};
  \node[emTag,above] at (1,1.97) {$#2$};
}
\newcommand{\EvenSideWeights}[4]{%
  \node[emTag,left] at (-1.53,1.25) {$#1$};
  \node[emTag,right] at (-.47,1.25) {$#2$};
  \node[emTag,left] at (.47,1.25) {$#3$};
  \node[emTag,right] at (1.53,1.25) {$#4$};
}
\newcommand{\EvenPictureBegin}{%
  \begin{tikzpicture}[x=8.8mm,y=8.8mm,baseline=(current bounding box.center)]
  \path[use as bounding box] (-2.1,-2.98) rectangle (2.1,2.98);
}
\newcommand{\EvenPictureEnd}{\end{tikzpicture}}

\newcommand{\EvenNoTunnelKnot}{%
  \EvenPictureBegin
  \EvenPlatStubs{emP}{emP}{emP}{emP}{emP}{emP}{emP}{emP}
  \EvenPlatBox\EvenTopCapLabels{p}{p}
  \node[emTag] at (0,-2.32) {$p=q$};
  \EvenPictureEnd}
\newcommand{\EvenNoTunnelLink}{%
  \EvenPictureBegin
  \EvenPlatStubs{emP}{emP}{emQ}{emQ}{emP}{emP}{emQ}{emQ}
  \EvenPlatBox\EvenTopCapLabels{p}{q}
  \node[emTag] at (0,-2.32) {$K_1\sqcup K_2$};
  \EvenPictureEnd}
\newcommand{\EvenOneTunnelKnot}{%
  \EvenPictureBegin
  \EvenPlatStubs{emP}{emP}{emQ}{emQ}{emP}{emQ}{emP}{emQ}
  \EvenPlatBox\EvenTopCapLabels{J_1:p}{J_2:q}\EvenLowerTunnel{d}
  \EvenPictureEnd}
\newcommand{\EvenOneTunnelLink}{%
  \EvenPictureBegin
  \EvenPlatStubs{emP}{emP}{emQ}{emQ}{emP}{emP}{emQ}{emQ}
  \EvenPlatBox\EvenTopCapLabels{K_1:p}{K_2:q}\EvenLowerTunnel{d}
  \EvenPictureEnd}
\newcommand{\EvenOneTunnelCoincident}{%
  \EvenPictureBegin
  \EvenPlatStubs{emP}{emP}{emEdge}{emEdge}{emP}{emEdge}{emP}{emEdge}
  \EvenPlatBox\EvenTopCapLabels{J_1:p}{J_2:2}\EvenLowerTunnel{m}
  \EvenPictureEnd}
\newcommand{\EvenTwoTunnelKnot}{%
  \EvenPictureBegin
  \EvenPlatStubs{emP}{emEdge}{emEdge}{emQ}{emP}{emEdge}{emEdge}{emQ}
  \EvenPlatBox\EvenPlatEdgeNames{e_1}{e_2}{e_4}{e_3}
  \EvenUpperTunnel{2}\EvenLowerTunnel{m}\EvenSideWeights{p}{2}{2}{q}
  \EvenPictureEnd}
\newcommand{\EvenTwoTunnelLinkOuter}{%
  \EvenPictureBegin
  \EvenPlatStubs{emP}{emEdge}{emEdge}{emQ}{emP}{emEdge}{emEdge}{emQ}
  \EvenPlatBox\EvenPlatEdgeNames{L_1}{L_2}{R_2}{R_1}
  \EvenUpperTunnel{2}\EvenLowerTunnel{m}\EvenSideWeights{p}{2}{2}{q}
  \EvenPictureEnd}
\newcommand{\EvenTwoTunnelLinkInner}{%
  \EvenPictureBegin
  \EvenPlatStubs{emP}{emEdge}{emQ}{emEdge}{emP}{emEdge}{emQ}{emEdge}
  \EvenPlatBox\EvenPlatEdgeNames{L_1}{L_2}{R_2}{R_1}
  \EvenUpperTunnel{2}\EvenLowerTunnel{m}\EvenSideWeights{p}{2}{q}{2}
  \EvenPictureEnd}
\newcommand{\EvenOrderTwoKnot}{%
  \EvenPictureBegin
  \EvenPlatStubs{emQ}{emChain}{emEdge}{emChain}{emQ}{emEdge}{emChain}{emChain}
  \EvenPlatBox\EvenPlatEdgeNames{e_1}{e_2}{e_4}{e_3}
  \EvenUpperTunnel{2}\EvenLowerTunnel{m}\EvenSideWeights{q}{2}{2}{2}
  \EvenPictureEnd}
\newcommand{\EvenOrderTwoLink}{%
  \EvenPictureBegin
  \EvenPlatStubs{emChain}{emChain}{emEdge}{emQ}{emChain}{emChain}{emEdge}{emQ}
  \EvenPlatBox\EvenPlatEdgeNames{L_1}{L_2}{R_2}{R_1}
  \EvenUpperTunnel{2}\EvenLowerTunnel{m}\EvenSideWeights{2}{2}{2}{q}
  \EvenPictureEnd}

\newcommand{\EvenLocalVertex}[1]{%
  \begin{tikzpicture}[x=10mm,y=10mm,baseline=(current bounding box.center)]
    \path[use as bounding box] (-1.35,-1.5) rectangle (1.35,1.65);
    \coordinate (O) at (0,0);
    \foreach \a/\n in {90/1,210/2,330/3}{
      \draw[emEdge] (\a:.66)--(\a:1.24);
      \node[emTag] at (\a:1.45) {$l_{\n}$};}
    \ifnum#1=0
      \draw[emGhost] (0,0) circle (.66);
      \draw[emGhost] (-.66,0) arc[start angle=180,end angle=360,x radius=.66,y radius=.18];
      \draw[draw=emGrey,line width=.55pt] (-.66,0) arc[start angle=180,end angle=0,x radius=.66,y radius=.18];
      \foreach \a in {90,210,330}{\draw[emEdge] (O)--(\a:.66);}
      \fill[emInk] (O) circle (1.8pt);
    \else
      \path[emCut] (0,0) circle (.66);
      \draw[emGhost] (-.66,0) arc[start angle=180,end angle=0,x radius=.66,y radius=.18];
      \draw[draw=emGrey,line width=.55pt] (-.66,0) arc[start angle=180,end angle=360,x radius=.66,y radius=.18];
      \foreach \a in {90,210,330}{
        \draw[emGhost] (O)--(\a:.66);
        \fill[emInk] (\a:.66) circle (1.55pt);}
      \node[emVertex] at (O) {};
    \fi
  \end{tikzpicture}}

\newcommand{\EvenTriangleCapped}{%
  \begin{tikzpicture}[x=9mm,y=9mm,baseline=(current bounding box.center)]
    \path[use as bounding box] (-2.5,-1.65) rectangle (2.5,1.7);
    \draw[emP] (-1.9,0)..controls(-1.5,1.62)and(1.5,1.62)..(1.9,0);
    \draw[emEdge] (-1.9,0)--(1.9,0);
    \draw[emQ] (-1.9,0)..controls(-1.5,-1.62)and(1.5,-1.62)..(1.9,0);
    \node[emVertex] at (-1.9,0) {};\node[emVertex] at (1.9,0) {};
    \node[emTag,above] at (0,1.22) {$p$};
    \node[emTag,above,fill=white] at (0,0) {$d$};
    \node[emTag,below] at (0,-1.22) {$q$};
  \end{tikzpicture}}
\newcommand{\EvenTriangleTruncated}{%
  \begin{tikzpicture}[x=9mm,y=9mm,baseline=(current bounding box.center)]
    \path[use as bounding box] (-2.7,-1.65) rectangle (2.7,1.7);
    \draw[emP] (-1.54,.56)..controls(-.9,1.44)and(.9,1.44)..(1.54,.56);
    \draw[emEdge] (-1.34,0)--(1.34,0);
    \draw[emQ] (-1.54,-.56)..controls(-.9,-1.44)and(.9,-1.44)..(1.54,-.56);
    \foreach \x in {-1.9,1.9}{
      \path[emCut] (\x,0) ellipse[x radius=.56,y radius=.78];
      \draw[emGhost] (\x,-.78) arc[start angle=-90,end angle=90,x radius=.18,y radius=.78];
      \draw[draw=emGrey,line width=.5pt] (\x,-.78) arc[start angle=-90,end angle=-270,x radius=.18,y radius=.78];
      \node[emVertex] at (\x,0) {};}
    \foreach \x/\y/\cx in {-1.54/.56/-1.9,-1.34/0/-1.9,-1.54/-.56/-1.9,1.54/.56/1.9,1.34/0/1.9,1.54/-.56/1.9}{
      \draw[emGhost] (\cx,0)--(\x,\y);
      \fill[emInk] (\x,\y) circle (1.55pt);}
    \node[emTag,above] at (0,1.22) {$p$};
    \node[emTag,above,fill=white] at (0,0) {$d$};
    \node[emTag,below] at (0,-1.22) {$q$};
  \end{tikzpicture}}

\newcommand{\EvenPanel}[3]{%
  \begin{minipage}[t]{#1}\centering
    #2\par\vspace{1mm}{\small #3}
  \end{minipage}}

\newcommand{\EvenFigureNoTunnels}{%
  \begin{figure}[htbp]\centering
    \EvenPanel{.43\linewidth}{\EvenNoTunnelKnot}{(a) A knot; $p=q$}
    \hfill
    \EvenPanel{.43\linewidth}{\EvenNoTunnelLink}{(b) Two components}
    \caption{Vertex-free models. Every component is marked. The box $B_r$
    contains a chosen four-plat braid for $K(r)$, with its full crossing
    data; it does not denote four trivial strands. The endpoint permutation
    is normalized to $(23)$ in (a) and to the identity in (b).
    All labels are even.}\label{fig:even-no-tunnels}
  \end{figure}}

\newcommand{\EvenFigureOneTunnel}{%
  \begin{figure}[htbp]\centering
    \EvenPanel{.32\linewidth}{\EvenOneTunnelKnot}{(a) Knot; distinct marks}
    \hfill
    \EvenPanel{.32\linewidth}{\EvenOneTunnelLink}{(b) Two components}
    \hfill
    \EvenPanel{.32\linewidth}{\EvenOneTunnelCoincident}{(c) Knot; coincident marks}
    \caption{One-tunnel models. In (a), $e_a=J_1$ and $e_b=J_2$;
    in (b), the marked edges are the two components. The surviving tunnel
    has weight $d\ge2$. In (c), $e_a=e_b=J_1$, $p=q$, and the tunnel
    has odd weight $m\ge3$; the other edge has weight $2$.
    Each arc $J_i$ includes its whole upper cap and both adjoining braid
    strands. Hollow vertices belong to the fully capped graph and are
    interpreted according to their incident weights.}\label{fig:even-one-tunnel}
  \end{figure}}

\newcommand{\EvenFigureTwoTunnels}{%
  \begin{figure}[htbp]\centering
    \EvenPanel{.32\linewidth}{\EvenTwoTunnelKnot}{(a) Knot; opposite edges}
    \hfill
    \EvenPanel{.32\linewidth}{\EvenTwoTunnelLinkOuter}{(b) Link; marks $L_1,R_1$}
    \hfill
    \EvenPanel{.32\linewidth}{\EvenTwoTunnelLinkInner}{(c) Link; marks $L_1,R_2$}
    \caption{Two-tunnel models, with weights $2$ and $m\ge3$ odd.
    In (a), $e_1,e_2,e_3,e_4$ occur cyclically on the knot,
    $\partial\tau_+=(e_1\cap e_2)\cup(e_3\cap e_4)$ and
    $\partial\tau_-=(e_2\cap e_3)\cup(e_4\cap e_1)$.
    The marked edges are $e_1,e_3$, of weights $p,q$; the other two
    have weight $2$. In (b)--(c), $K_1=L_1\cup L_2$ and
    $K_2=R_1\cup R_2$; the two possible placements on the second
    component are displayed. In every panel the box retains the
    spatial braid data.}\label{fig:even-two-tunnels}
  \end{figure}}

\newcommand{\EvenFigureOrderTwo}{%
  \begin{figure}[htbp]\centering
    \EvenPanel{.44\linewidth}{\EvenOrderTwoKnot}{(a) Knot; chain $e_2\cup e_3$}
    \hfill
    \EvenPanel{.44\linewidth}{\EvenOrderTwoLink}{(b) Link; chain $K_1$}
    \caption{The order-two axis chain when $p=2$, $q\ge4$, and both
    tunnels survive. The thick teal arcs are the entire image of
    $\operatorname{Ax}(a)$ in the capped graph; all have weight $2$.
    The marked $b$-edge has weight $q$. In (a), $e_b=e_1$ and the
    chain is $e_2\cup e_3$, although the representative $e_a=e_3$
    may be chosen opposite $e_b$. In (b), the chain is the whole
    component disjoint from $e_b$. The tunnel weights are $2$ and
    $m\ge3$ odd.}\label{fig:even-order-two-chain}
  \end{figure}}

\newcommand{\EvenFigureVertices}{%
  \begin{figure}[htbp]\centering
    \EvenPanel{.31\linewidth}{\EvenLocalVertex{0}}{(a) $s(v)>1$; retained vertex}
    \hfill
    \EvenPanel{.31\linewidth}{\EvenLocalVertex{1}}{(b) Euclidean piece ($s(v)=1$)}
    \hfill
    \EvenPanel{.31\linewidth}{\EvenLocalVertex{2}}{(c) Hyperbolic piece ($s(v)<1$)}
    \caption{Realization of a capped trivalent vertex with incident
weights $l_i\in\{2,3,\ldots,\infty\}$, where
$s(v)=l_1^{-1}+l_2^{-1}+l_3^{-1}$ and $1/\infty=0$.
In (a), the vertex is retained and the dotted sphere is its link.
In (b)--(c), a neighbourhood of the vertex is removed.
The three black marks record the incident germs. A finite label
gives a cone point of that order, while an infinite label gives
a puncture, represented by a boundary circle in the compact core.
With finite labels the resulting boundary is a closed Euclidean
turnover in (b) or a closed hyperbolic turnover in (c).
The infinite-label case in (b) is $(2,2,\infty)$ and gives
$D^2(2,2)$, a rank-one cusp section. The dashed radial segments
and hollow centre in (b)--(c) belong to the removed cap.
The drawings record incidence and truncation.}
    \label{fig:even-vertex-realization}
  \end{figure}}

\newcommand{\EvenFigureTriangle}{%
  \begin{figure}[htbp]\centering
    \EvenPanel{.46\linewidth}{\EvenTriangleCapped}{(a) Fully capped theta graph}
    \hfill
    \EvenPanel{.50\linewidth}{\EvenTriangleTruncated}{(b) Both formal vertices reopened}
    \caption{The fully capped and truncated models of
$\Delta(p,q,d)$ for integers $p,q,d\ge2$ satisfying
$p^{-1}+q^{-1}+d^{-1}<1$. The theta graph has marked edge
weights $p,q$ and tunnel weight $d$. Reopening its two vertices
gives a compact core with boundary turnovers $S^2(p,q,d)$.
The complete quotient is orbifold homeomorphic to
$S^2(p,q,d)\times\mathbb R$ and has infinite volume.
When $p,q$ are even, marked order-two surgery gives
$\Odih(0/1;1,d)$.}\label{fig:even-triangle-model}
  \end{figure}}

\hypersetup{pdftitle={Kleinian groups generated by two elliptics of even order}}
\theoremstyle{plain}
\newtheorem{theorem}{Theorem}[section]
\newtheorem{proposition}[theorem]{Proposition}
\newtheorem{lemma}[theorem]{Lemma}
\newtheorem{remark}[theorem]{Remark}
\newtheorem{corollary}[theorem]{Corollary}
\theoremstyle{definition}
\newtheorem{definition}[theorem]{Definition}
\newtheorem{hypothesis}[theorem]{Hypothesis}
\newtheorem{example}[theorem]{Example}

\newcommand{\HH}{\mathbb H}

\newcommand{\orb}{\mathrm{orb}}
\newcommand{\PSL}{\mathrm{PSL}_2(\mathbb C)}
\newcommand{\nc}[1]{\langle\!\langle #1\rangle\!\rangle}
\newcommand{\Odih}{\mathcal O}

\newcommand{\ax}{\operatorname{Ax}}

\DeclareMathOperator{\tr}{tr}
\title[Two elliptics of even order]{Kleinian groups generated by two elliptics of even order}
\author{Gaven Martin and Jeroen Schillewaert}
\date{15 November 2026}
\address[GM]{Institute for Advanced Study, Massey University, Auckland, New Zealand}
\address[JS]{Department of Mathematics, The University of Auckland, New Zealand}
\email[GM]{g.j.martin@massey.ac.nz}
\email[JS]{j.schillewaert@auckland.ac.nz}
\begin{document}

\begin{abstract}
Let $\Gamma<\PSL$ be a non-elementary discrete group generated by
two elliptics of even order, with a nontrivial relation in the
corresponding free product of cyclic groups. We classify the possible
compact orbifold cores of $\mathbb H^3/\Gamma$, together with their
cusp data and marked meridians. They are finite meridional fillings
of the two-bridge link and Heckoid models of
Akiyoshi--Ohshika--Parker--Sakuma--Yoshida. Their weighted graphs
have at most two standard tunnels, whose weights in the two-tunnel
case are $2$ and an odd integer. The original generators are powers
of primitive meridians which generate $\Gamma$, and we determine the
possible meridian placements for all even orders. The groups are
geometrically finite, their conformal boundary components are
turnovers, and each oriented core with its cusp data determines its
complete hyperbolic holonomy up to M\"obius conjugacy. The proof uses
trace rigidity, dihedral surgery, half-turn symmetries, and an
index-two cover.
\end{abstract}
\maketitle

\section{Introduction}\label{sec:introduction}

We consider non-elementary Kleinian groups generated by two elliptic
elements of even order, with a relation other than their order
relations. We determine their compact orbifold cores and the
positions of the generating meridians. The resulting models are
finite-weight versions of the two-bridge link and Heckoid models
in the classification of two-parabolic-generator groups by
Akiyoshi--Ohshika--Parker--Sakuma--Yoshida \cite{AOPSY}. We call
these finite-weight models \emph{elliptic Heckoid orbifolds},
including the meridional cone fillings of two-bridge links.

Write $C_n$ for the cyclic group of order $n$. For an ordered
generating pair $(g_1,g_2)$ of exact orders $p_0,q_0$, the
\emph{marked epimorphism} $C_{p_0}*C_{q_0}\to\Gamma$ sends the
distinguished factor generators to $g_1,g_2$. It is injective on
each factor, a property also called \emph{smoothness}. A primitive
elliptic meridian is a once-around rotation in the full pointwise
stabilizer of its axis in $\Gamma$. For an elliptic $g$, write
$\ax(g)$ for its fixed geodesic in $\mathbb H^3$.

\begin{theorem}[Main theorem]\label{thm:classification}
Let $\Gamma=\langle g_1,g_2\rangle<\PSL$ be discrete and
non-elementary. Suppose that $g_1,g_2$ have even orders $p_0,q_0$
and that the marked epimorphism $C_{p_0}*C_{q_0}\longrightarrow\Gamma$
is not injective. There are primitive elliptic meridians
$a,b\in\Gamma$ such that $\Gamma=\langle a,b\rangle,
 g_1\in\langle a\rangle, g_2\in\langle b\rangle,
  p_0\mid p, q_0\mid q$,
where $p,q$ are the orders of $a,b$.
The group $\Gamma$ is M\"obius conjugate to the complete hyperbolic
holonomy of an oriented compact orbifold $Q$ of one of the forms
\[
 \mathcal L(r)_{p,q},\qquad
 \mathcal M_0(r;d)_{p,q},\qquad
 \mathcal M_1(r;m)_p,\qquad
 \mathcal M_2(r;m)_{p,q}.
\]
Here $r\in\mathbb Q$, $d\ge2$ is an integer and $m\ge3$ is odd.
The four families are defined in Section~\ref{sec:heckoid} and
listed in Table~\ref{tab:group-models}. In particular,
$\Gamma\cong\pi_1^{\orb}(Q)$, with $a,b$ represented by the
specified marked meridians. The possible placements are given in
Theorem~\ref{thm:mainC}.

The group is geometrically finite, and every component of its
conformal boundary is a sphere with three cone points. The oriented
compact orbifold $Q$, together with the boundary components
designated as cusps, determines this holonomy up to M\"obius conjugacy.
\end{theorem}

\begin{remark}
The same conclusion holds when $g_1,g_2$ have arbitrary finite
orders, provided the full pointwise stabilizer of each of their
axes in $\Gamma$ has even order. Their primitive meridians then
have the even orders required by the spatial classification.
Without this proviso the conclusion can fail when a generator has
odd order, see Proposition \ref{prop:odd-counterexamples}.
\end{remark}

Here $K(r)$ denotes the two-bridge knot or link of rational slope
$r$, in the convention of Section~\ref{sec:weighted-models}.
The notation $\mathcal M_i$ follows \cite[Definition~3.4]{AOPSY}.
The subscripts record the finite orders replacing its infinite
meridional weights. The family $\mathcal L$ consists of the
two-bridge link fillings. Section~\ref{sec:heckoid} defines these
four families and gives the correspondence with the diagrams
and the usual Heckoid index in \cite{AOPSY}.

\begin{table}[ht]
\caption{The families in Theorem~\ref{thm:classification}, with
even marked orders $p,q$. Only models admitting the stated complete
non-elementary hyperbolic realization occur. The definitions are in
Section~\ref{sec:heckoid}, and the figures show the marked edges.}
\label{tab:group-models}
\centering\small
\begin{tabular}{@{}p{.26\textwidth}p{.43\textwidth}p{.23\textwidth}@{}}
\toprule
Family & Parameters & Figure\\
\midrule
$\mathcal L(r)_{p,q}$ &
$r\in\mathbb Q$, with $p=q$ if $K(r)$ is a knot &
\ref{fig:even-no-tunnels}\\[4pt]
$\mathcal M_0(r;d)_{p,q}$ &
$r\in\mathbb Q$, integer $d\ge2$ &
\ref{fig:even-one-tunnel}(a)--(b)\\[4pt]
$\mathcal M_1(r;m)_p$ &
$K(r)$ a knot, odd $m\ge3$, and $q=p$ &
\ref{fig:even-one-tunnel}(c)\\[4pt]
$\mathcal M_2(r;m)_{p,q}$ &
$r\in\mathbb Q$, odd $m\ge3$ &
\ref{fig:even-two-tunnels}\\
\bottomrule
\end{tabular}
\end{table}

For integers $p,q,d\ge2$ with $1/p+1/q+1/d<1$, the planar
theta graph with weights $p,q,d$ gives the triangle group
\[
 \pi_1^{\orb}\bigl(\mathcal M_0(0;d)_{p,q}\bigr)
 =\Delta(p,q,d)
 =\langle a,b\mid a^p=b^q=(ab)^d=1\rangle.
\]
Its quotient is homeomorphic to $S^2(p,q,d)\times\mathbb R$ and
has infinite volume. Example~\ref{rem:triangle-example} explains
this realization. Example~\ref{ex:two-tunnel-lattice} gives a finite-volume
two-tunnel model with marked orders $(4,6)$. Thus both the integral
and the strictly half-integral Heckoid families have even-order
elliptic realizations.

The proof has three steps. First,
Proposition~\ref{thm:rigid} and
Lemma~\ref{lem:rigid-boundary-turnover} give geometric finiteness
and turnover conformal boundary. Lowering the two marked weights
to two then permits the dihedral surgery and classification of
\cite{AOPSY}. This identifies an embedded two-bridge graph with
standard tunnels. Second, half-turn symmetries and an index-two
cover determine the possible marked edges. The additional
arguments for $(p,q)=(4,4)$ and $(2,4)$ use a midpoint-stabilizer
obstruction and the lifted spatial graph, respectively. Finally,
restoring the marked weights recovers the compact orbifold with
its cusp data. Turnover rigidity, in
Proposition~\ref{prop:pared-rigidity}, determines its complete
holonomy up to M\"obius conjugacy.

The vertex labels determine the ends and finite volume
(Corollary~\ref{cor:ends}). The resulting cusp list gives a
commensurability obstruction and, together with Hoffman's covering
restrictions, excludes finite covers by hyperbolic knot complements
in $S^3$ (Corollaries~\ref{cor:cusp-consequences}
and~\ref{cor:no-knot-cover}).

Let $O=\mathbb H^3/\Gamma$, and write $|O|$ for its underlying
topological space. A \emph{fully capped graph} is obtained by
embedding the compact core in $S^3$, filling its underlying
boundary spheres by balls, and extending the singular arcs inside
these balls. The prescribed deletion of cap vertices and infinite
edges recovers the compact core with its cusp data, as explained in
Section~\ref{sec:weighted-models}. The following consequences
describe the vertex-free and disconnected cases.

\begin{theorem}[Vertex-free case]\label{thm:mainA}
Let $\Gamma=\langle a,b\rangle$ satisfy the hypotheses of
Theorem~\ref{thm:classification}, with $a,b$ primitive of orders $p,q$.
Suppose its fully capped graph $\Pi$ has no vertex. Then $\Gamma$
is a cocompact lattice, $|O|=S^3$, and the singular locus is a knot
or a two-component link $L\subset S^3$. Every component is marked,
and hence every label is even. The group
$\Gamma=\pi_1^\orb(S^3,L,\mathbf n)$ is generated by the two marked
meridians, where $\mathbf n$ lists the component labels.
Moreover, $L$ is a two-bridge knot or link.
\end{theorem}

\begin{corollary}\label{cor:disconnected-even}
Under the same group hypotheses, if the fully capped graph is
disconnected, it is a two-bridge two-component link, with one mark
on each component. In particular each component is unknotted.
\end{corollary}

Earlier classifications impose different restrictions.
Klimenko--Kopteva \cite{KK} assume real traces for the two
generators and their commutator, and no invariant plane.
Maclachlan--Martin \cite{MaclachlanMartin} assume that the
commutator of the elliptic generators is elliptic.
Martin--Salehi--Yamashita \cite{MSY} classify the arithmetic
lattices with generator orders $(4,p)$. Here the parity of the
generating orders gives the dihedral quotient used in the proof.
Section~\ref{sec:parity-limits} exhibits compact hyperbolic examples
with generating orders $(3,4)$ and $(3,5)$ whose tetrahedral
singular graphs lie outside the four families above.

Sections~\ref{sec:setup} and~\ref{sec:weighted-models} establish
the rigidity and orbifold preliminaries. The spatial classification
is stated in Section~\ref{sec:mainresults} and proved in
Sections~\ref{sec:warmup}--\ref{sec:ordertwo}.
Section~\ref{sec:rigidity} recovers the holonomy and gives the
examples and consequences for ends and covers.

{\bf Acknowledgements:} Both authors are funded by the New Zealand Marsden Fund. The majority of the research was carried out whilst we enjoyed being hosted by the Max Planck Institute in Bonn in August 2026 jointly with Eric Chesebro.
His help was crucial in understanding the methods used in \cite{AOPSY}. 
The main ideas are all due to the authors, the write-up of this paper was aided by OpenAI's Astra 6.0 and Anthropic's Fable 5.2. All AI assisted content was
revised by the authors.

\section{Trace rigidity and conformal boundary}\label{sec:setup}

\subsection{Standing conventions}

All orbifolds are effective, orientable, and without reflector boundary.
At an interior point, a $3$-orbifold is a ball modulo a finite orientation-preserving
rotation group, acting faithfully. A point of a singular edge has cyclic
local group; a spherical vertex has a noncyclic finite rotation group.
Their fundamental groups are orbifold fundamental groups unless the
underlying topological space is explicitly indicated. Thus local
rotations are retained in the group, although their loops are killed
on passage to the underlying space. A Kleinian group is a discrete
subgroup of $\PSL$; it is non-elementary if its limit set has more than
two points. The limit set is the accumulation set on
$\partial\HH^3\cong\widehat{\mathbb C}$ of an orbit in $\HH^3$. Write $C_n$ for the
cyclic group of order $n$, and
\[
 D_n=\langle u,v\mid u^2=v^2=(uv)^n=1\rangle
 \quad(n\ge1),\qquad D_\infty=C_2*C_2.
\]
Thus $|D_n|=2n$, $D_1=C_2$, and $D_2=C_2\times C_2$.

A meridian is the image of a once-around local loop about a singular
edge, with either orientation and an arbitrary path to the base point.
In a hyperbolic orbifold, a primitive elliptic meridian of order $n$
has rotation angle $\pm2\pi/n$. 

For general orbifold and hyperbolic geometry background, see
\cite[Chapter~13]{Thurston} and \cite{Ratcliffe}.

\begin{hypothesis}\label{hyp:main}
The group $\Gamma=\langle a,b\rangle<\PSL$ is discrete and
non-elementary. The elements $a,b$ are primitive
elliptic meridians of orders $p,q$, and the marked epimorphism
$C_p*C_q\twoheadrightarrow\Gamma$ has nontrivial kernel.
\end{hypothesis}

Put $O=\mathbb H^3/\Gamma$. We write $(M_0,P)$ for a compact
orbifold core $M_0$ in which $P\subset\partial M_0$ records the cusp
sections. Its existence follows from
\cite[Theorems~5.1--5.2]{AOPSY}. Once geometric finiteness is established, the components of
$\partial M_0\setminus P$, with their boundary circles adjacent to
$P$ replaced by punctures, identify with the components of
$\Omega(\Gamma)/\Gamma$. We use these identifications to give the
nonparabolic boundary its conformal structures.
Section~\ref{sec:weighted-models} gives the compactification
and realization conventions in full.

\subsection{Density and marked rigidity}\label{sec:density}

The domain of discontinuity is
$\Omega(\Gamma)=\widehat{\mathbb C}\setminus\Lambda(\Gamma)$, and
$\Omega(\Gamma)/\Gamma$ is the conformal boundary. We use the usual
equivalent definitions of geometric finiteness; one is that a fixed
positive neighbourhood of the convex core has finite volume. The
convex core is the quotient of the hyperbolic convex hull of the limit
set. A marked representation retains the identification with the
abstract group, and algebraic convergence means convergence of the
images of each group element. A marked quasiconformal deformation is
a conjugacy on the Riemann sphere by a quasiconformal homeomorphism,
with M\"obius conjugate representations identified; quasiconformality
allows uniformly bounded infinitesimal distortion \cite{Maskit}.

\begin{theorem}[Namazi--Souto, {\cite[Theorem~12.1, pp.~387--388]{NS}}]
\label{orbifold-density}
Let $\Gamma<\PSL$ be discrete, finitely generated and
non-elementary, with torsion allowed. There are discrete faithful representations
$\rho_n:\Gamma\to\PSL$ with geometrically finite images such that
$\rho_n(\gamma)\to\gamma$ for every $\gamma\in\Gamma$.
\end{theorem}

\begin{proof}[Deduction from the cited theorem]
In \cite[\S12.1]{NS}, the space $AH(\Gamma)$ is defined for groups
with torsion and consists of conjugacy classes of discrete faithful
representations of the fixed group $\Gamma$. Their proof chooses
a finite-index torsion-free normal subgroup $\Gamma_0$ and
geometrically finite approximations to its inclusion whose ending
invariants are equivariant under $G=\Gamma/\Gamma_0$. The same
proof extends these approximations to discrete faithful
representations $\rho_n$ of $\Gamma$. Their images are geometrically
finite because they are finite extensions of the geometrically
finite groups $\rho_n(\Gamma_0)$.

The convergence on $\Gamma_0$ determines convergence on $\Gamma$.
Choose loxodromics $u,v\in\Gamma_0$ with disjoint fixed-point sets.
For large $n$, the attracting and repelling fixed points of
$\rho_n(u)$ and the attracting fixed point of $\rho_n(v)$ form a
triple converging to the corresponding three distinct points for
$u,v$. For $\alpha\in\Gamma$, normality gives
$\alpha u\alpha^{-1},\alpha v\alpha^{-1}\in\Gamma_0$.
The conjugacy identities force $\rho_n(\alpha)$ to send the first
triple to the corresponding triple for
$\rho_n(\alpha u\alpha^{-1}),\rho_n(\alpha v\alpha^{-1})$.
The latter converges to the image of the first limiting triple
under $\alpha$. A M\"obius transformation depends continuously on
its values at three distinct points, so $\rho_n(\alpha)\to\alpha$.
\end{proof}

\paragraph{Trace coordinates.}\label{ex:characters}
For a representation $\rho:\Gamma\to\PSL$, choose lifts
$A,B\in\mathrm{SL}_2(\mathbb C)$ of $\rho(a),\rho(b)$ and put
\[
x=\tr A,\qquad y=\tr B,\qquad z=\tr(AB).
\]
The lifts define a representation of the free group on $a,b$.
A relation of $\Gamma$ evaluates on them to $I$ or $-I$.
Changing the signs of $A,B$ changes the trace triple to
$(\epsilon x,\eta y,\epsilon\eta z)$, where
$\epsilon,\eta\in\{1,-1\}$.
For fixed finite projective orders of $A,B$, the possible values of
$x,y$ form finite sets. Thus $x,y$ are constant on every connected
family on which the lifts are chosen continuously. The sign changes
above describe the finite ambiguity in choosing those lifts.
The remaining coordinate $z$ is constrained by the mixed relation
in the following proposition.
Every word trace is an integer polynomial
in $x,y,z$, and an irreducible pair is determined up to simultaneous
conjugacy by this triple \cite{Goldman2009,Horowitz}. Here
irreducible means that the two matrices have no common invariant
line in $\mathbb C^2$.

\begin{proposition}\label{thm:rigid}\label{lem:pinned}
Let $\Gamma=\langle a,b\rangle<\PSL$ be discrete and non-elementary,
where $a,b$ are nontrivial elliptics of orders $p,q$. If the marked
epimorphism $C_p*C_q\to\Gamma$ has nontrivial kernel, then the
character of the inclusion is isolated in the projective character
variety of $\Gamma$. Moreover, $\Gamma$ is geometrically finite and
has no nontrivial marked quasiconformal deformation.
\end{proposition}
\begin{proof}
Choose lifts $A,B$ of $a,b$ and write
$x=\tr A$, $y=\tr B$, $z=\tr(AB)$. The pair is irreducible because
$\Gamma$ is non-elementary.
The kernel of $C_p*C_q\to\Gamma$ meets each factor trivially.
Consequently a nonidentity kernel element has a cyclically reduced
conjugate of the form
\[
R=a^{i_1}b^{j_1}\cdots a^{i_k}b^{j_k},
\qquad 1\le i_\nu<p,\quad 1\le j_\nu<q,\quad k\ge1.
\]
Here cyclic reduction is performed in the free product, combining
end syllables that belong to the same factor until the first and
last syllables belong to different factors. With $x,y$ fixed,
the word trace is a polynomial $P_R(Z)\in\mathbb Q(x,y)[Z]$.
We verify its degree after this specialization.
\label{ex:relation-degree}
Let $\lambda^{\pm1},\mu^{\pm1}$ be the eigenvalues of $A,B$.
Put $S_0(X)=0$, $S_1(X)=1$ and
$S_{j+1}(X)=XS_j(X)-S_{j-1}(X)$ for $j\ge1$.
Cayley--Hamilton gives $C^2-(\tr C)C+I=0$ for
$C\in\mathrm{SL}_2(\mathbb C)$, and induction yields
\[
 C^j=S_j(\tr C)C-S_{j-1}(\tr C)I\qquad(j\ge1).
\]
Evaluating this recurrence on the two eigenvalues gives
\[
S_j(x)=\frac{\lambda^j-\lambda^{-j}}{\lambda-\lambda^{-1}},
\qquad
S_j(y)=\frac{\mu^j-\mu^{-j}}{\mu-\mu^{-1}}.
\]
Evaluate the word on
\[
A_\lambda=\begin{pmatrix}\lambda&1\\0&\lambda^{-1}\end{pmatrix},
\qquad
B_\mu(T)=\begin{pmatrix}\mu&0\\T&\mu^{-1}\end{pmatrix}.
\]
Their product trace is $T+c$, with
$c=\lambda\mu+(\lambda\mu)^{-1}$. Each matrix
$A_\lambda^iB_\mu(T)^j$ is affine in $T$, with coefficient
\[
S_j(y)\begin{pmatrix}S_i(x)&0\\\lambda^{-i}&0\end{pmatrix}.
\]
Thus the coefficient of $T^k$ in the word trace, and hence the
leading coefficient of $P_R$, is
\begin{equation}\label{eq:relation-leading-coefficient}
\prod_{\nu=1}^{k}S_{i_\nu}(x)S_{j_\nu}(y)
=
\prod_{\nu=1}^{k}
\frac{\lambda^{i_\nu}-\lambda^{-i_\nu}}{\lambda-\lambda^{-1}}
\frac{\mu^{j_\nu}-\mu^{-j_\nu}}{\mu-\mu^{-1}}
\ne0.
\end{equation}
Indeed, $\lambda^2,\mu^2$ have exact orders $p,q$, and all the
exponents lie strictly between zero and the corresponding order.
Therefore $\deg P_R=k$.
For every representation of $\Gamma$ with these lifted generator
traces, the relation $R=1$ implies $P_R(z)=2$ or $P_R(z)=-2$.
Only this necessary consequence of the matrix relation is used.
There are finitely many possible values of $z$.
In a sufficiently small neighbourhood of the inclusion, continuous
lifts have constant traces $x,y$, since every representation satisfies
$a^p=b^q=1$ and these traces lie in finite sets. The neighbouring
pairs remain irreducible, and the finite set of possible product
traces isolates the inclusion character.
By Theorem~\ref{orbifold-density}, there are discrete faithful
geometrically finite representations $\phi_n:\Gamma\to\PSL$
converging algebraically to the inclusion. Choose generator lifts
$A_n\to A$, $B_n\to B$. For all sufficiently large $n$ their traces
equal $x,y$, while $z_n=\tr(A_nB_n)$ belongs to the finite zero set
of $P_R^2-4$ and tends to $z$. Hence $z_n=z$ eventually. Equality
of the irreducible characters gives simultaneous conjugacy, so
$\Gamma$ is geometrically finite.
Finally, a marked quasiconformal deformation is joined to the
inclusion by scaling its invariant Beltrami coefficient. Choose
lifts continuously along this path. Their traces $x,y$ are constant,
and the product trace takes values in the same finite set.
The character is constant, so every representation on the path
is M\"obius conjugate to the inclusion.
\end{proof}

\begin{lemma}\label{lem:rigid-boundary-turnover}
Under the hypotheses of Proposition~\ref{thm:rigid}, every
component of the conformal boundary is a hyperbolic turnover,
with punctures allowed.
\end{lemma}

\begin{proof}
The group is geometrically finite by Proposition~\ref{thm:rigid}.
Let $\partial_0M_0$ be the nonparabolic boundary of the compact
pair $(M_0,P)$, with boundary circles adjacent to $P$ replaced
by punctures.
Its components carry the conformal structures of
$\Omega(\Gamma)/\Gamma$ and are hyperbolic $2$-orbifolds of finite
type by Ahlfors' finiteness theorem \cite{Maskit}. If this boundary
is empty the conclusion is immediate.
For a component $F$, let $\mathcal T(F)$ be its space of marked
conformal structures, with cone orders and punctures fixed,
modulo conformal equivalence respecting the marking up to isotopy.
The orbifold Ahlfors--Bers theorem
\cite[Definition~8.4.2 and Proposition~8.4.4]{BP}
realizes small paths in $\prod_F\mathcal T(F)$ by continuous
quasiconformal conjugacies $f_t$, with $f_0=1$.
Suppose such a path starts at the given conformal boundary.
Proposition~\ref{thm:rigid} makes its holonomy character constant.
The conjugating M\"obius maps $m_t$ may be chosen continuously
near $t=0$, with $m_0=1$. To see this, diagonalize the first
generator using a continuous eigenbasis and normalize the
upper-right entry of the second generator to $1$ by a continuous
diagonal conjugation. The distinct eigenvalues and irreducibility
permit these choices. The two trace equations then determine the
diagonal entries, and the determinant determines the remaining
off-diagonal entry. Thus the normalized pair is constant. The maps $m_t^{-1}f_t$ commute with
$\Gamma$ and descend to an isotopy of $\Omega(\Gamma)/\Gamma$.
Since $m_t$ conformally identifies the quotient for the deformed
holonomy with the original quotient, the boundary point is
constant in $\prod_F\mathcal T(F)$. Every factor therefore has
dimension zero.
For a hyperbolic component of genus $g$ with $n$ cone points or
punctures, the complex dimension is $3g-3+n$. Dimension zero
forces $g=0,n=3$, since the alternative $g=1,n=0$ is Euclidean.
Thus each component is a sphere with three cone points or
punctures. Its lack of conformal moduli can also be seen by
sending the three specified points to $0,1,\infty$.
\end{proof}

\section{Orbifold cores, weighted graphs and surgery}
\label{sec:weighted-models}

\paragraph{Local groups and their links.}\label{ex:local}
Let $x$ be an interior point of an orientable $3$-orbifold. A sufficiently
small neighbourhood is $B^3/H_x$, where $H_x<\mathrm{SO}(3)$ is finite.
Its boundary is the \emph{link} $S^2/H_x$ of $x$.
A regular point has $H_x=1$.
If $H_x=C_n$, the quotient ball contains a singular diameter of order
$n$, and its link is a sphere with two cone points of order $n$.
Write $S^2(l_1,\ldots,l_k)$ for the sphere with the indicated
cone-point orders, interpreting an infinite order as a puncture.
At a noncyclic spherical vertex the possibilities are
\[
\begin{array}{c|c|c}
H_x&|H_x|&S^2/H_x\\ \hline
D_n\ (n\ge2)&2n&S^2(2,2,n)\\
A_4&12&S^2(2,3,3)\\
S_4&24&S^2(2,3,4)\\
A_5&60&S^2(2,3,5).
\end{array}
\]
These are the rotation groups of a regular polygon, tetrahedron,
octahedron and icosahedron. The three cone points record the three
orbits of directions with nontrivial stabilizer. Thus an edge germ at
a vertex is an orbit of rays, not an individual axis in the covering
ball. This distinction will determine when a projected axis folds.

\subsection{Weighted models and their realization}

We now describe the compact pair $(M_0,P)$ used above.
A cusp neighbourhood is a quotient of a horoball by its
parabolic stabilizer. \emph{Cusp truncation} removes
open such neighbourhoods and retains their cross-sections. The
\emph{relative compactification} also adds boundary at the remaining
ends, compactifying the cusp cross-sections at the same time. It gives
a compact pair $(M_0,P)$: $M_0$ is an orbifold compact core, and
$P\subset\partial M_0$ records precisely the cusp sections. Its interior
recovers $O$, and the core inclusion induces an isomorphism on
fundamental groups; no meridional filling relation is introduced by
truncation. Existence and the pared properties follow from
\cite[Theorems 5.1--5.2]{AOPSY}.

Here \emph{pared} has the meaning of \cite[Definition~3.1]{AOPSY}:
$M_0$ is compact, irreducible and has a finite manifold cover;
$P$ is a union of incompressible annular and toric $2$-orbifolds;
rank-two free abelian subgroups are peripheral; and embedded annular
orbifolds with essential boundary in $P$ are parallel into $P$.
Incompressible means that inclusion is injective on orbifold fundamental
groups; peripheral means conjugate into the group of a component of $P$.
Irreducibility means that every embedded spherical $2$-suborbifold bounds
a discal $3$-orbifold, a quotient of a ball by a finite rotation group.
The singular locus of $O$ is distinguished from the weighted
graph used to record this compactification.

A weighted graph specifies a compact orbifold by the local orders
on its finite edges and by the neighbourhoods to be removed.
We call the graph a \emph{model}, and call the resulting compact
orbifold together with its specified parabolic locus its
\emph{realization}. A \emph{formal vertex} is a vertex whose
neighbourhood is removed in this realization. We use
\cite[Conventions~6.1--6.2]{AOPSY}. A weighted
graph $(W,\Pi,w)$ consists of a compact orientable $3$-manifold, a finite
properly embedded graph, and weights in
$\{2,3,\ldots\}\cup\{\infty\}$. Circle components are allowed.
Boundary vertices are univalent. Interior vertices are trivalent,
with the additional permitted case of a quadrivalent vertex all
of whose incident weights are $2$.
At a trivalent vertex with incident weights $l_1,l_2,l_3$, put
\[
 s(v)=\frac1{l_1}+\frac1{l_2}+\frac1{l_3},\qquad 1/\infty=0.
\]
The vertex is spherical, Euclidean, or hyperbolic according as $s(v)>1$,
$s(v)=1$, or $s(v)<1$. To form the actual compact orbifold, remove
neighbourhoods of the infinite-weight edges and of the nonspherical
vertices, and retain each finite edge with its indicated local order.

\EvenFigureVertices

A quadrivalent all-$2$ vertex is removed and records a pillow cusp
$S^2(2,2,2,2)$, the sphere with four cone points of order two.
A cusp has rank one or two according to the rank of the translation
subgroup of its parabolic stabilizer. Compact rank-one sections are
annuli or $D^2(2,2)$, a disc with two order-two cone points.
The latter has a twofold annulus cover and orbifold Euler
characteristic zero. Compact rank-two sections are tori, pillows,
or closed Euclidean turnovers. The parabolic locus is specified
as in the cited conventions.

For a sphere component $S$ of $\partial W$, condition {\rm(SC)} means that
its induced punctured cone-sphere is neither bad, spherical, discal, nor an
annulus; equivalently, $S$ has at least three endpoints of $\Pi$ and, if it
has exactly three, with incident weights $l_1,l_2,l_3$, then
\[
  \frac1{l_1}+\frac1{l_2}+\frac1{l_3}\le 1,
  \qquad 1/\infty=0 .
\]
In particular the boundary pattern $(2,2,\infty)$ is allowed.
Its compact realization is $D^2(2,2)$, whose interior is the
punctured orbifold $S^2(2,2,\infty)$. This is the condition in
\cite[Convention~6.1]{AOPSY}. A bad orientable closed
$2$-orbifold is a sphere with one cone point, or with two cone points of
unequal orders; it has no surface cover. Condition {\rm(SC)} is not a
substitute for the requirement, imposed separately, that the source pair
in an orbifold surgery be pared.

Throughout the classification statements, $\Pi$ denotes a
\emph{fully capped model}. All boundary spheres of the underlying
compact core have been filled by balls in the ambient manifold $W$.
Turnover pieces and pillow cusps are represented by formal vertices,
and parabolic annular data is represented by infinite-weight edges.
Thus $\partial W=\varnothing$. The realization recovers the compact
pair $(M_0,P)$ by removing the prescribed neighbourhoods, so its
underlying space may have boundary. Such models exist by
Lemma~\ref{lem:fullycapped}. All classification assertions concern
any chosen fully capped model of the specified pared core with its
prescribed parabolic data.

Choose open singular edges of $O$ contained in the images of
$\ax(a)$ and $\ax(b)$, and denote the corresponding capped edges
by $e_a,e_b$. Choose orientations and meridian base paths so that
their meridians represent $a,b$. The projected axis of an order-two
meridian can follow several successive singular edges, and the
choice of $e_a$ or $e_b$ selects one of these edge pieces.
If the chosen edges coincide, then $p=q$.

Weight changes and augmentation are performed in the sense of
\cite[Definition 6.3]{AOPSY}. In the fully capped convention,
{\rm(SC)} is vacuous on $\partial W$. Equivalently, the operations may
be performed in the uncapped presentation, where newly spherical
turnover boundaries are capped by discal orbifolds.

\begin{lemma}[Meridional presentation]\label{lem:labelfree}
Let $U=\pi_1(W\setminus\Pi)$ and choose one meridian $x_e$ for each
edge. For a weighted realization $Q$ as above,
\[
 \pi_1^{\orb}(Q)=
 U\big/\nc{x_e^{w(e)}:w(e)<\infty}.
\]
Crossing and vertex relations belong to $U$ and do not depend on the
weights. Capping a spherical turnover boundary by its discal orbifold
does not change the fundamental group.
\end{lemma}

\begin{proof}
Filling a retained edge neighbourhood adds precisely the relation
$x_e^{w(e)}=1$. There is no power relation on a drilled edge. The local
vertex relations already hold in the complement. At a capped spherical
turnover, its fundamental group maps isomorphically to that of the
discal cap, so orbifold van Kampen adds no further relation.
\end{proof}

\begin{lemma}\label{lem:background}
Let $\Gamma<\PSL$ be a finitely generated Kleinian group generated
by elliptic elements, and let $(M_0,P)$ be a compact core obtained by relative
compactification of $\mathbb H^3/\Gamma$. Then $|M_0|$ is simply
connected and every component of its boundary is a sphere.
Capping these spheres by balls gives $S^3$. In particular, every
fully capped weighted-graph realization of this core has closed
ambient manifold $S^3$.
\end{lemma}

\begin{proof}
Put $O=\mathbb H^3/\Gamma$ and $Y=|M_0|$.
The quotient of $S^2$ by a finite subgroup of $\mathrm{SO}(3)$
has underlying space $S^2$, so the underlying space of each
interior chart $B^3/H$ is a ball. Boundary collars give half-ball
charts. Hence $Y$ is a compact orientable $3$-manifold.
Armstrong's theorem \cite{Armstrong} gives
\[
 \pi_1(|O|)=\Gamma/\nc{\text{elliptic elements}}=1,
\]
because $\Gamma$ is generated by elliptic elements.
Relative compactification and cusp truncation preserve this
underlying fundamental group, so $Y$ is simply connected.
Poincar\'e--Lefschetz duality
\cite[Theorem~3.43, p.~254]{Hatcher} and the exact sequence
of the pair give
\[
 H_2(Y,\partial Y;\mathbb Q)\longrightarrow
 H_1(\partial Y;\mathbb Q)\longrightarrow H_1(Y;\mathbb Q).
\]
The first group is $H^1(Y;\mathbb Q)=0$ and the last is zero.
Thus $H_1(\partial Y;\mathbb Q)=0$, and every boundary component
is a sphere. Capping these spheres by balls preserves simple
connectivity. The Poincar\'e theorem identifies the resulting
closed manifold with $S^3$
\cite{KleinerLott,Perelman1,Perelman2,Perelman3}.
\end{proof}

\begin{lemma}[Fully capped models]\label{lem:fullycapped}
Under Hypothesis~\ref{hyp:main}, the pared compactification $(M_0,P)$
admits a fully capped weighted-graph model $(S^3,\Pi,w)$ in the sense
of \cite[Conventions 6.1--6.2]{AOPSY}.
\end{lemma}
\begin{proof}
The pair $(M_0,P)$ is pared by \cite[Theorems~5.1--5.2]{AOPSY}.
By Lemma~\ref{lem:background}, every component of its underlying
boundary is a sphere. By Lemma~\ref{lem:rigid-boundary-turnover},
each component of $\overline{\partial M_0\setminus P}$, with its boundary
circles replaced by punctures, is a hyperbolic turnover.
The components of $P$ are annuli, discs $D^2(2,2)$, Euclidean
turnovers, or pillows. A torus component cannot lie on an
underlying boundary sphere.
Let $S$ be a boundary component containing a rank-one cusp section.
Cut $|S|$ along the boundary circles of these sections. The dual
graph of the resulting pieces is a tree, since disjoint circles
on a sphere separate and cutting along each circle increases
the number of pieces by one. Suppress the bivalent vertices
corresponding to annular components of $P$. The remaining vertices
correspond to the nonparabolic pieces and the $D^2(2,2)$ components.
A nonparabolic piece with $c$ cone points has $3-c$ boundary
circles. A $D^2(2,2)$ component has two cone points and one boundary
circle. Choose a collar of this boundary circle within each
$D^2(2,2)$ component.
Thicken this tree to a ball $B_S$, using one ball for each vertex
and one connecting $1$-handle for each edge. Its boundary pieces
and connecting annuli have the same incidence pattern as those
on $S$. Choose their boundary identifications to respect the cone
points and the original gluing circles. They give an identification
$\partial B_S\cong |S|$ respecting this entire boundary pattern.
For an edge ending at a $D^2(2,2)$ vertex, identify the boundary
annulus of its handle with the chosen collar in that component.
The thickening is a ball because its underlying graph is a tree.
In each vertex ball, join its centre to the cone points and to
the centres of its attaching discs. Continue the latter arcs
along the handle cores. Give every handle-core edge weight
$\infty$, and give an arc ending at a cone point its cone order.
Every new vertex is trivalent. At a nonparabolic piece its link
is the corresponding hyperbolic turnover, with punctures recorded
by the infinite edges. At a $D^2(2,2)$ piece its type is
$(2,2,\infty)$. Attach the finite arcs to the corresponding
singular arcs of $M_0$. A boundary component containing no rank-one cusp section is a
closed hyperbolic turnover, a Euclidean turnover, or a pillow.
Cap it by a ball and cone its cone points to one vertex.
Use a trivalent vertex for a turnover and a quadrivalent all-$2$
vertex for a pillow.
Let $\Pi$ be the union of the singular locus of $M_0$ and the
graphs just added, with the original weights retained on $M_0$.
The resulting closed ambient manifold is $S^3$ by
Lemma~\ref{lem:background}. The new vertices are nonspherical,
and each capping ball is a regular neighbourhood of the new
vertices and infinite edges that it contains. Deleting these
neighbourhoods therefore recovers $M_0$, with its finite singular
arcs, up to boundary collars. The meridional annuli of the infinite
edges recover the annular cusp sections. At a terminal
$(2,2,\infty)$ vertex, the Euclidean disc piece together with
its adjacent annular collar recovers the prescribed $D^2(2,2)$
section. The closed Euclidean turnovers and pillows recover the
remaining components of $P$. Thus the realization recovers the
specified pared pair $(M_0,P)$.
\end{proof}

\subsection{Dihedral orbifolds and surgery inputs}

We use the oriented rational-tangle convention of \cite[Section~2]{AOPSY}:
$K(r)$ is obtained by gluing the tangles of slopes $\infty$ and $r$
along their common marked four-punctured sphere. If $r=u/v$ is in
lowest terms with $v>0$, then $K(r)$ is a knot when $v$ is odd and a
two-component link when $v$ is even. The core tunnels of the two
tangles are denoted by $\tau_+$ and $\tau_-$, respectively. These
conventions specify the spatial embedding and its orientation.
Thus the denominator denoted by $p$ in \cite{AOPSY} is $v$ here.
We reserve $p,q$ for the marked elliptic orders.

We use the standard Klein-four group action on the two rational
tangles, glued along their marked boundary sphere
\cite[\S2 and Figure~4]{AOPSY}. Let $\alpha_0$ be the
slope-zero separating loop on the marked boundary of the upper
tangle. The horizontal involution fixes the upper tunnel
pointwise. The other two nonidentity elements exchange its
endpoints. The vertical involution is the one preserving the
orientation of $\alpha_0$. The product of the horizontal and
vertical involutions is the planar involution.

For coprime positive integers $d_+,d_-$, write
\[
 \Odih(r;d_+,d_-)=
 (S^3,K(r)\cup\tau_+\cup\tau_-,w),\qquad
 w|_{K(r)}=2,\quad w(\tau_\pm)=d_\pm.
\]
A rational tangle is a pair of simultaneously boundary-parallel arcs in a ball with
a specified marking of its four boundary endpoints; the slope records
this marking relative to the trivial tangle. A standard tunnel joins
the two arcs in the corresponding bridge ball. A tunnel of weight $1$ is deleted and resulting bivalent vertices are
suppressed. This is a spherical orbifold. Interchanging the tunnels may
change the slope; we allow this change and normalize $d_+\le d_-$.

\begin{theorem}[Akiyoshi--Ohshika--Parker--Sakuma--Yoshida]
\label{thm:AOPSYinput}
The following results hold.
\begin{enumerate}[(a)]
\item Surgery on a pared orbifold in the sense of
\cite[Definition 6.3]{AOPSY} produces a very good compact orbifold with no bad $2$-suborbifold and
no spherical boundary component \cite[Lemma 6.4]{AOPSY}.
\item A compact orientable $3$-orbifold with nonempty singular set,
no bad $2$-suborbifold, no spherical boundary component, and dihedral
fundamental group is one of $\Odih(r;d_+,d_-)$, the all-$2$ two-component
unlink in $S^3$, the weight-$2$ unknot in $\mathbb{RP}^3$, or
$D^2(2,2)\times I$ \cite[Theorem 4.1]{AOPSY}.
\item Let $d_+,d_-\ge2$ be coprime.
No orientation-preserving involution of $\Odih(r;d_+,d_-)$
preserves a non-tunnel edge setwise
\cite[Corollary 12.7(2)]{AOPSY}. The model has the standard Klein-four
symmetries, transitive on the four link edges
\cite[Proposition 12.6 and Figure 14]{AOPSY}.
\end{enumerate}
\end{theorem}

An orbifold is \emph{good} if it has a manifold cover, and
\emph{very good} if it has a finite-sheeted manifold cover.
Goodness is equivalent to injectivity of every local group
into the orbifold fundamental group. Indeed, the local groups
of the universal orbifold cover are the kernels of these maps.
They are all trivial exactly when that universal cover is a
manifold. A manifold cover has a manifold universal cover,
which proves the equivalence.

\begin{definition}\label{def:ordertwo}
The marked order-two surgery replaces the weights on $e_a,e_b$
by $2$ and takes the weighted realization, including the spherical
vertices created by the change. Denote the resulting actual compact
orbifold by $O_2$.
\end{definition}
This is an orbifold surgery in the sense of
\cite[Definition~6.3]{AOPSY} for arbitrary finite marked orders.
The source realization is pared, the closed ambient manifold
makes {\rm(SC)} vacuous, and the weights incident to a pillow
vertex remain $2$.

To describe augmentation in the uncapped core, consider a
trivalent vertex whose original neighbourhood was removed.
Its three new labels may have reciprocal sum greater than one.
The corresponding boundary is then a spherical turnover
$S^2/H$, and augmentation attaches the discal orbifold $B^3/H$
along it. In the fully capped graph this same operation retains
the vertex neighbourhood that the new weights make spherical.
When $p,q$ are even, Lemma~\ref{lem:paritysurgery} identifies
$\pi_1^{\orb}(O_2)$ with $\Gamma/\nc{a^2,b^2}$.
The divisibilities $2\mid p,q$ enter in this quotient formula.

\subsection{Elliptic Heckoid models}\label{sec:heckoid}

We use the notation $\mathcal M_i$ of
\cite[Definition~3.4]{AOPSY}, adding a subscript for the finite
orders that replace its infinite meridional weights.
Throughout this subsection $r\in\mathbb Q$, $p,q\ge2$,
$d\ge2$ is an integer, and $m\ge3$ is odd. The following are
weighted graphs in $S^3$, realized by the convention in
Section~\ref{sec:weighted-models}.

\begin{enumerate}[(i)]
\item $\mathcal L(r)_{p,q}$ is the link $K(r)$ with weights
$p,q$ on its two components. If $K(r)$ is a knot, it has the
single weight $p=q$. For a two-component link, choose one marked
meridian on each component. For a knot, both marked meridians are
supported on its single component.

\item $\mathcal M_0(r;d)_{p,q}$ has graph $K(r)\cup\tau_-$,
with weight $d$ on the tunnel. The closures of the components of
$K(r)\setminus\partial\tau_-$ are its two non-tunnel edges.
Give these edges weights $p,q$ and mark one meridian on each.
For a knot they are two arcs with the same endpoints. For a
two-component link each is a loop based at one tunnel endpoint.

\item $\mathcal M_1(r;m)_p$ is defined when $K(r)$ is a knot.
Its graph is $K(r)\cup\tau_-$, with weights $p,2$ on the two
knot arcs and weight $m$ on the tunnel. Both marked meridians
are supported on the weight-$p$ edge, so $q=p$.

\item $\mathcal M_2(r;m)_{p,q}$ has graph
$K(r)\cup\tau_+\cup\tau_-$, with tunnel weights $2,m$ in that
order. Two non-tunnel edges with no common endpoint carry the
marked weights $p,q$, and the other two have weight $2$.
For a knot the marked edges are opposite arcs. For a two-component
link they lie on different components, and the standard pair is
the pair preserved setwise by the vertical involution, as in
\cite[Definition~3.4(3)]{AOPSY}.
\end{enumerate}

For reference in the spatial proof, label the four arcs of a knot
in the two-tunnel case cyclically as $e_1,e_2,e_3,e_4$, with
\[
 \partial\tau_+=(e_1\cap e_2)\cup(e_3\cap e_4),\qquad
 \partial\tau_-=(e_2\cap e_3)\cup(e_4\cap e_1).
\]
The standard marked pair is $e_1,e_3$.
For a two-component link write $K(r)=K_1\sqcup K_2$ and
$K_j=J_{1,j}\cup J_{2,j}$, with the vertical involution
interchanging $J_{i,1}$ and $J_{i,2}$ (the indexing of \cite[Section 9, Subcase 3.1]{AOPSY}). After applying a horizontal
involution, a disjoint marked pair is $J_{1,1},J_{i,2}$ for
$i=1$ or $2$. The second choice gives the standard model with
slope $r+1$, by \cite[Remark~3.5(3)]{AOPSY}. The relevant
homeomorphism preserves orientation and the ordered tunnels
\cite[Proposition~2.1(1a)]{AOPSY}. The individual marked weights
are carried with their edges.

These four families are the \emph{elliptic Heckoid models}.
Their compact orbifolds are obtained by retaining spherical
vertices and removing neighbourhoods of nonspherical vertices.
The Euclidean boundary components are designated as cusps.
The classification concerns the models admitting a complete
non-elementary hyperbolic realization with these cusp data.

The associated groups are determined directly by the spatial
graphs. For a model $Q$ represented by $(S^3,\Lambda,w)$,
Lemma~\ref{lem:labelfree} gives
\begin{equation}\label{eq:heckoid-graph-group}
 \pi_1^{\orb}(Q)=
 \pi_1(S^3\setminus\Lambda)\big/
 \nc{\mu_e^{\,w(e)}:e\in E(\Lambda)},
\end{equation}
where $\mu_e$ is a meridian of the edge $e$.
Thus a diagram gives a presentation using its crossing and vertex
relations and the displayed meridional power relations.
In the one-tunnel family this has the familiar form
\[
 \pi_1^{\orb}\bigl(\mathcal M_0(r;d)_{p,q}\bigr)
 =\langle f,g\mid f^p=g^q=w_r(f,g)^d=1\rangle.
\]
Here $f,g$ are the meridians of the upper trivial tangle, and
$w_r$ is the boundary of a compressing disc for the lower
rational tangle, read in this free meridian basis. The presentation
follows from van Kampen as in \cite[Section~2]{LeeSakuma}.
Setting $d=1$ gives the corresponding presentation for
$\mathcal L(r)_{p,q}$, with $p=q$ in the knot case.

\paragraph{The published models and our figures.}
Replace the infinite weights in \cite[Figure~1]{AOPSY} by the
marked finite orders. Its panels~(1), (2), (3), and (4) then give,
respectively, Figure~\ref{fig:even-no-tunnels},
Figure~\ref{fig:even-one-tunnel}(a)--(b),
Figure~\ref{fig:even-one-tunnel}(c), and
Figure~\ref{fig:even-two-tunnels}. The edges of weight $2$ and
the tunnel weights are retained. The vertex realization is then
recomputed from the incident weights, as in
Figure~\ref{fig:even-vertex-realization}. The two link placements
in Figure~\ref{fig:even-two-tunnels}(b)--(c) correspond to
\cite[Figure~5(3)--(4)]{AOPSY} and to the slopes $r,r+1$ just
described. Figure~\ref{fig:even-order-two-chain} additionally
records the projected axis when a marked order is two.

For comparison with the traditional Heckoid index, write
$S(s;n)$ for the parabolic Heckoid orbifold of
\cite[Definition~3.2]{AOPSY} and
\cite[Definition~3.2]{LeeSakuma}. Omitting the finite-order
subscript from $\mathcal M_i$ denotes its parabolic model.
The description in \cite[Propositions~5.2--5.3]{LeeSakuma},
in the form preceding \cite[Definition~3.4]{AOPSY}, gives
\begin{equation}\label{eq:heckoid-parameter-dictionary}
\begin{aligned}
 \mathcal M_0(r;d)&\cong S(r;d),\\
 \mathcal M_1(u/v;m)&\cong S(2u/v;m/2)
 &&\text{if }v\text{ is odd},\\
 \mathcal M_2(u/v;m)&\cong S(u/(2v);m/2)
 &&\text{if }u\text{ is odd}.
\end{aligned}
\end{equation}
Here $u/v$ is reduced with $v>0$, and the last line uses the
standard marked pair. For a knot with even numerator, the change
$r\mapsto r+1$ makes the numerator odd and preserves the ordered
tunnels. The opposite pair is determined up to the standard
symmetry \cite[Remark~3.5(2)]{AOPSY}. For a two-component link
the reduced numerator is already odd. The displayed formulas
follow by substituting $s=2u/v$ and $s=u/(2v)$ into the published
forward conversion. The Heckoid orbifold $S(s;n)$ is called even
when $n$ is integral and odd when $n$ is strictly half-integral.
This terminology is independent of the parity of the finite
marked orders.

\section{Spatial classification}\label{sec:mainresults}

We now classify the fully capped embedded graph and the positions of
the two marked meridians. The order-two surgery determines the ambient
two-bridge tunnel model; the subsequent sections determine which
marked weights and placements can be restored. The end data is part
of this identification, so that Section~\ref{sec:rigidity} can recover
the complete hyperbolic holonomy.

\begin{theorem}[Spatial classification]\label{thm:mainC}
Let $\Gamma=\langle a,b\rangle<\PSL$ be discrete and non-elementary,
where $a,b$ are primitive elliptic meridians of even orders $p,q$.
Suppose that the marked epimorphism
$C_p*C_q\twoheadrightarrow\Gamma$ is not injective. Then $\Gamma$
is geometrically finite and has no rank-one or pillow cusps.
Its pared compact core admits a fully capped weighted-graph model
with closed ambient manifold $S^3$. For any chosen fully capped
model $\Pi$ with the prescribed parabolic data, the following
assertions hold.

\begin{enumerate}[(I)]
\item There exist $r\in\mathbb Q$ and coprime integers
$d_+,d_-\ge1$ such that
\[
O_2\cong\Odih(r;d_+,d_-).
\]
The original pared orbifold $(M_0,P)$ is recovered by restoring the
weights $p,q$ on the marked edges and reopening precisely those capped
vertices that become Euclidean or hyperbolic after the restoration.
After changing the identification of the exceptional unknotted
all-$2$ theta graph by its threefold symmetry, if necessary, the
marked edges lie on $K(r)$ and not on the tunnels. If $\Pi$ has a
vertex, $d_+=d_-=1$ cannot occur.

\item Suppose that $p,q\ge4$ and $\Pi$ has a vertex.
After normalizing $1\le d_+\le d_-$:
\begin{enumerate}[(a)]
\item Let $d_+=1<d_-$. The model is
$\mathcal M_0(r;d_-)_{p,q}$ with distinct marked edges, or
$K(r)$ is a knot, $d_-$ is odd, $p=q$, and the model is
$\mathcal M_1(r;d_-)_p$ with coincident marked edges.

\EvenFigureOneTunnel

\item If $d_+\ne1$, then $d_+=2$ and $d_-=m\ge3$ is
odd. The model is $\mathcal M_2(r;m)_{p,q}$, allowing the slope
change described in Section~\ref{sec:heckoid}.
For a two-component link, in the notation of that subsection,
the marked edges can be written
\[
 e_a=J_{1,1},\qquad e_b=J_{i,2}\quad(i\in\{1,2\}).
\]
For a knot, use the cyclic edge labelling of that subsection.
The marked edges are opposite and may be labelled
\[
 e_a=e_1,\qquad e_b=e_3.
\]
\end{enumerate}
In particular, a two-tunnel model with both tunnel weights odd does not
occur.

\EvenFigureTwoTunnels

\item If $\Pi$ has no vertex, then $\Gamma$ is cocompact and
$O=\HH^3/\Gamma$ is a labelled two-bridge knot or link in $S^3$.
Every link component is marked. If $\Pi$ is disconnected, it has
no vertex and is a two-component two-bridge link, with both
components unknotted and one mark on each component.

\EvenFigureNoTunnels

\item If one marked order is two, interchange $a,b$ if necessary
so that $p=2$; then $q\ge4$. The one-tunnel conclusion in
\textup{(II)(a)} holds. If both tunnels survive, normalize
$d_+\le d_-$; then
$d_+=2$ and $d_-=m\ge3$ is odd. If $K(r)$ has two components,
the image of $\ax a$ traverses the entire component not containing
$e_b$. If $K(r)$ is a knot, choose the cyclic labelling in
\textup{(II)(b)} with $e_b=e_1$. Then the image of $\ax a$
traverses precisely $e_2\cup e_3$, and the marked edge for $a$
may be chosen to be the opposite edge $e_3$.

\EvenFigureOrderTwo

\end{enumerate}
\end{theorem}

\paragraph{Drawing conventions.}
The figures represent fully capped weighted graphs in $S^3$.
Hollow symbols denote vertices of these graphs. Their incident
weights determine whether the vertices are retained or reopened,
as in Figure~\ref{fig:even-vertex-realization}.
Colour and thickness distinguish marked edges, and numerical
labels give their local orders. A tunnel of weight $1$ is absent.

The box $B_r$ contains the braid of a four-plat presentation of
the oriented tunnel pair of $K(r)$, including all its crossings.
Number its endpoints $1,2,3,4$ from left to right at each end.
The upper and lower bridge pairings are $(12)(34)$.
Cap-preserving endpoint motions normalize the braid permutation
to $(23)$ for a knot and to the identity for a two-component link.
These motions extend over the bridge balls and preserve their
tunnels. The slope $r$ refers to the oriented rational-tangle
convention of Section~\ref{sec:weighted-models}.

\section{Even-labelled links and mod-two homology}\label{sec:warmup}

For a link with even labels, lowering the labels to two gives
a dihedral quotient that recognizes the two-bridge locus.
For the weighted graphs, mod-two homology constrains the
supporting edges of the marked generators.

\subsection{The even-labelled link theorem}

For a link $L=K_1\cup\cdots\cup K_c\subset S^3$ with meridians
$\mu_1,\ldots,\mu_c$, and integers $n_i\ge2$, write
$\mathbf n=(n_1,\ldots,n_c)$ and put
\[
 G(L;\mathbf n)=\pi_1(S^3\setminus L)/\nc{\mu_i^{n_i}:1\le i\le c}.
\]
Alexander duality \cite[Theorem~3.44, p.~255]{Hatcher} gives
$H_1(S^3\setminus L;\mathbb Z)\cong\mathbb Z^c$, with basis the
meridians, equivalently the linking-number coordinates. Imposing the
power relations therefore gives abelianization $\bigoplus_i C_{n_i}$,
with the meridian of $K_i$ in the $i$-th summand. In particular that meridian has exact order $n_i$.
Write $G_2(L)=G(L;2,\ldots,2)$.

\begin{lemma}\label{lem:paritydescent}\label{lem:dihedral}\label{lem:sigma2}
If every $n_i$ is even, there is a meridian-preserving epimorphism
$G(L;\mathbf n)\twoheadrightarrow G_2(L)$. If $G(L;\mathbf n)$ is generated
by two meridians $a,b$, write $\bar a,\bar b$ for their images
in $G_2(L)$. These images are involutions, and
\[
 G_2(L)\cong D_n\quad(1\le n\le\infty),\qquad
 \pi_1(\Sigma_2(L))\cong\langle\bar a\bar b\rangle,
\]
where $\Sigma_2(L)$ is the double branched cover and the displayed
rotation subgroup has index two.
\end{lemma}

\begin{proof}
The inclusion $\nc{\mu_i^{n_i}}\subseteq\nc{\mu_i^2}$ gives the
epimorphism. The total parity map sends every meridian of $G_2(L)$
to $1\in C_2$. The images $\bar a,\bar b$ are consequently nontrivial
involutions. Every even word in them is a power of $\bar a\bar b$,
so the parity kernel is cyclic and the group is dihedral. The cover
defined by parity is a manifold, since each local group $C_2$ maps
injectively to the deck group; it is exactly $\Sigma_2(L)$.
\end{proof}

\begin{theorem}\label{thm:warmup}
Let $L\subset S^3$ be a nonempty non-split link, with component
labels $n_i\ge2$. If $G(L;\mathbf n)$ is generated by two
meridians of even order, then every component label is even,
$L$ has at most two components, and $L$ is two-bridge, with the
unknot allowed. Its double branched cover is a lens space or $S^3$.
\end{theorem}
\begin{proof}
The abelianization $\bigoplus_i C_{n_i}$ is generated by the two
marked meridian classes. Every component therefore carries a mark,
and its label is the even order of that meridian.
By the preceding lemma, $(S^3,L;2,\ldots,2)$ has dihedral fundamental
group. It contains no bad $2$-suborbifold: an embedded sphere meets a
link in an even number of points, and all its cone points have equal
order $2$. It is closed. Apply Theorem~\ref{thm:AOPSYinput}(b).
The weight-two unknot in $\mathbb{RP}^3$ has a different underlying
closed manifold. The product $D^2(2,2)\times I$ has nonempty boundary,
whereas the present orbifold is closed. The two-component unlink is
split and is excluded by hypothesis. In a spherical model, a surviving tunnel produces
trivalent vertices. Hence both tunnel weights are $1$, and the link is
$K(r)$. The double cover of each rational tangle branched over its
two arcs is a solid torus. Gluing these two solid tori gives its standard
lens-space cover (including $S^3$); the $S^2\times S^1$ gluing is the
two-component unlink case already excluded. The component bound also follows immediately from
$H_1(G_2(L);\mathbb F_2)=\mathbb F_2^c$, generated by two classes.
\end{proof}

Locally along the branch link, the underlying double cover is
$(z,t)\mapsto(z^2,t)$ on $D^2\times I$. The all-$2$ labelling
makes this an orbifold cover with trivial upstairs edge stabilizer,
which explains why the parity kernel is a manifold group.

\subsection{Mod-2 homology}

\begin{lemma}\label{lem:mod2weighted}
Let $Q$ be the realization of a weighted graph $(S^3,\Pi,w)$.
Write
$H_1^{\orb}(Q;\mathbb F_2)
=(\pi_1^{\orb}Q)^{\mathrm{ab}}\otimes_{\mathbb Z}\mathbb F_2$.
Let $\Pi_{\rm ev}$ be the spanning subgraph obtained by deleting
all edges of finite odd weight. Thus it retains every edge of
even or infinite weight. Then
\[
H_1^{\orb}(Q;\mathbb F_2)
\cong H_1(S^3\setminus\Pi_{\rm ev};\mathbb F_2)
\cong H^1(\Pi_{\rm ev};\mathbb F_2).
\]
Under these isomorphisms, the class of a retained edge meridian
corresponds to the cohomology class of its edge cochain.
A meridian of an edge of finite odd weight is zero. At a trivalent
vertex with one finite odd germ and two even or infinite germs,
the latter two meridians have the same mod-two homology class.
\end{lemma}
\begin{proof}
The finite-weight case is \cite[Lemma~6.5]{AOPSY}.
For the stated extension, abelianise the presentation of
Lemma~\ref{lem:labelfree} over $\mathbb F_2$.
A finite odd power relation kills its meridian,
an even power relation imposes no condition, and an infinite edge
has no power relation. Alexander duality identifies
$H_1(S^3\setminus\Pi;\mathbb F_2)$ with
$H^1(\Pi;\mathbb F_2)$, sending a meridian to the class of the
cochain supported on its edge. Killing the odd-edge cochains
leaves the cochain presentation of $H^1(\Pi_{\rm ev};\mathbb F_2)$.
Applying Alexander duality to $\Pi_{\rm ev}$ proves both
isomorphisms \cite[Corollary~8.7, p.~353]{Bredon}.
At a trivalent vertex the sum of the three incident meridian
classes is zero. When one germ has finite odd weight, its class
vanishes and the other two classes are equal. A cut edge of
$\Pi_{\rm ev}$ has zero class, since its edge cochain is the
coboundary of the indicator of one component of its complement.
\end{proof}

\section{The vertex-free and disconnected cases}\label{sec:vertexfree}

\begin{lemma}\label{lem:vertexfree}
If a labelled link group in $S^3$ is generated by two meridians, each
component carries a marked meridian and there are at most two components.
For one component the marked orders coincide and the meridians are
conjugate up to inversion; for two components they lie on distinct
components.
\end{lemma}

\begin{proof}
In the abelianization $\bigoplus_i C_{n_i}$, a meridian has nonzero
coordinate only in its own component summand. Each summand must
therefore be represented by a marked meridian. The remaining assertions
follow from the definition of a meridian.
\end{proof}

\begin{proof}[Proof of Theorem~\ref{thm:mainA}]
There is no turnover end or pillow cusp, since these are recorded by
vertices in a fully capped model. The only possible remaining boundary
components of its realization would be tori from infinite-weight circle
components. Lemma~\ref{lem:background} excludes such boundary
components. Thus the realization is closed, $\Gamma$ is cocompact,
and $|O|=S^3$.

Its singular set is a labelled link. By Lemma~\ref{lem:vertexfree}
it has at most two components, each carrying a marked meridian, so all
labels are even. A hyperbolic orbifold is irreducible. A splitting sphere
for this link would have to bound a nonsingular ball on one side, which
is incompatible with a splitting into two nonempty sublinks. Hence the
link is non-split, and Theorem~\ref{thm:warmup} applies.
\end{proof}

\section{The connected case with vertices}\label{sec:connected}

Throughout this section assume Hypothesis~\ref{hyp:main},
with $p,q$ even. Connectedness and the presence of vertices are imposed
only where stated.

\subsection{Surgery and vertex tracking}

\begin{lemma}\label{lem:paritysurgery}
Assume Hypothesis~\ref{hyp:main}, with $p,q$ even. Then
\[
 \pi_1^{\orb}(O_2)=\Gamma/\nc{a^2,b^2}
\]
is a nontrivial dihedral group.
\end{lemma}

\begin{proof}
Evenness of $p,q$ and Lemma~\ref{lem:labelfree} give the quotient
formula. The quotient is generated by two elements of order at most
two. The construction in Definition~\ref{def:ordertwo} is an
orbifold surgery. By Theorem~\ref{thm:AOPSYinput}(a), $O_2$
is very good. Its marked edges give a nonempty singular
set, so local-group injectivity rules out the trivial group.
Every nontrivial quotient of $C_2*C_2$
is $D_n$ for some $1\le n\le\infty$.
\end{proof}

\begin{lemma}\label{lem:vertex-tracking}
For any fully capped weighted model, lowering finite weights
to integers at least two preserves every trivalent vertex
of the capped graph. If its new reciprocal sum exceeds $1$, it is an
interior singular vertex of the new realization; otherwise it remains
nonspherical boundary data. If an incident edge has infinite weight,
this data is a punctured turnover piece and need not be an entire
boundary component. A quadrivalent all-$2$ vertex persists
as a pillow boundary component under the marked order-two surgery.
\end{lemma}

\begin{proof}
The underlying graph and its incidences are unchanged. A lowered weight
can only increase the reciprocal sum. In the fully capped convention,
changing from a nonspherical to a spherical vertex retains its cone;
in the uncapped convention this is precisely augmentation. A pillow
vertex already has all weights $2$, so no incident weight changes.
\end{proof}

\begin{lemma}\label{lem:degenerate}
Assume Hypothesis~\ref{hyp:main}, with $p,q$ even. Then $O_2\cong\Odih(r;d_+,d_-)$ for some $r\in\mathbb Q$
and coprime positive $d_+,d_-$. The original fully capped graph has
no infinite-weight edge or pillow vertex. If it has a vertex, the
tunnel weights are not both $1$.
\end{lemma}

\begin{proof}
Apply Theorem~\ref{thm:AOPSYinput}(a),(b) and
Lemma~\ref{lem:paritysurgery}. If $O_2$ is closed, the original
model has no infinite edge or pillow vertex, since these data
are unchanged by surgery and would still give boundary.
Every trivalent vertex is then spherical after surgery and
is an actual trivalent vertex of $O_2$.
The weight-$2$ unknot in $\mathbb{RP}^3$ is excluded because the
closed ambient manifold is $S^3$, by Lemma~\ref{lem:background}.
If $O_2$ is the all-$2$ two-component unlink, the original graph
also has no vertices. Restoring the marked weights gives a
closed hyperbolic orbifold with a split singular link. Its
splitting sphere contradicts irreducibility.

Suppose that $O_2=D^2(2,2)\times I$. There is no interior
trivalent singular vertex, so every trivalent vertex of the
capped graph is nonspherical after surgery. Restoring the
original weights only decreases its reciprocal sum, so every
such vertex was nonspherical before surgery as well.
The two realizations therefore delete the same vertices,
infinite edges and pillow vertices. Their underlying ball
and two trivial singular arcs agree.
The images of $a,b$ generate
$H_1^{\orb}(O_2;\mathbb F_2)\cong\mathbb F_2^2$, so the two
marks lie on different arcs. Restoring the weights gives
$\Gamma\cong C_p*C_q$, with $a,b$ generating conjugates of the
two respective cyclic factors. Their fixed vertices in the Bass--Serre tree are distinct,
since they have the two different vertex types. Let $I$ join
these vertices. A nonidentity element of either factor sends
the first edge of $I$ at its fixed vertex to a different edge,
because edge stabilizers are trivial. The half-trees determined
by an interior point of an edge of $I$ give ping-pong for the two
factors \cite[Chapter~I, \S4]{Serre}. Thus the marked map
$C_p*C_q\to\Gamma$ is injective, contrary to
Hypothesis~\ref{hyp:main}.
The remaining model is $\Odih(r;d_+,d_-)$ and is spherical,
hence closed. The first observation excludes infinite edges
and pillow vertices in the original model. If $d_+=d_-=1$,
the resulting singular graph has no vertices, so the original
fully capped graph also has no vertices.
\end{proof}

\begin{proof}[Proof of Corollary~\ref{cor:disconnected-even}]
By Lemma~\ref{lem:degenerate}, the graph is the singular graph of a
spherical dihedral model. If a tunnel survives, that graph is connected.
Thus both tunnel weights are $1$, and the graph is the two-component
link $K(r)$. Its components are unknotted by its two-bridge position.
Lemma~\ref{lem:vertexfree} puts one mark on each component.
\end{proof}

\begin{corollary}\label{cor:reopen}
The original weighted realization is recovered from $O_2$ by restoring
the weights $p,q$ and truncating exactly those vertices whose restored
reciprocal sum is at most $1$.
\end{corollary}

\begin{proof}
Restoring weights decreases their reciprocals and exactly reverses the
construction in Lemma~\ref{lem:vertex-tracking}.
\end{proof}

\subsection{Excluding marked tunnels}

\begin{proposition}\label{prop:tunnelexclusion}
Suppose that, under an identification
\[
O_2\cong\Odih(r;d_+,d_-),
\]
one marked edge is a tunnel. Then the other tunnel is absent and
$K(r)$ is the unknot. Hence $O_2$ is the unknotted theta graph with
all three weights equal to $2$. Its symmetry permutes the three
edges, so the identification may be changed so that both marked edges
lie on the two non-tunnel edges.
\end{proposition}

\begin{proof}
Write $r=u/v$ in lowest terms with $v>0$. Let $\bar a,\bar b$
be the images of the marked generators in
$\pi_1^{\orb}(O_2)$, and suppose that $e_a$ is a tunnel.
Its weight is $2$. Denote the weight of the other tunnel by
$d$, where $d=1$ means that this other tunnel is absent.

Kill $\bar a$. Van Kampen identifies this quotient with the
result of deleting the tunnel $e_a$ and concatenating the
two link germs at each of its endpoints. Locally the vertex
relation identifies the two remaining meridians up to inverse
and conjugacy. Since both have order-two relations, the
concatenated edge again has weight $2$. Thus the resulting
weighted graph is $\Odih(r;1,d)$ or $\Odih(r;d,1)$, according
to which tunnel was erased. By
\cite[Proposition~12.2, equation~(5)]{AOPSY},
\[
\pi_1^{\orb}(O_2)/\nc{\bar a}\cong D_{vd}.
\]
The same quotient is generated by the image of $\bar b$,
whose order is at most two. Since $|D_{vd}|=2vd$, it follows
that $v=d=1$.
Consequently $K(r)$ is the unknot, the other tunnel is absent,
and $O_2$ is the unknotted theta graph with all three weights
equal to $2$. The quotient just computed is $D_1=C_2$.
If both marks had been on $e_a$, then both generators would
have been killed, giving the trivial quotient. Hence the
marked edges are distinct. The threefold orientation-preserving
symmetry of the unknotted theta graph sends the unmarked edge
to the tunnel position, placing both marked edges on $K(r)$.
\end{proof}

Proposition~\ref{thm:rigid} and Lemma~\ref{lem:fullycapped} give the
geometric finiteness and model assertions of Theorem~\ref{thm:mainC};
Lemma~\ref{lem:degenerate} excludes rank-one and pillow cusps.
Part (I) now follows from Lemma~\ref{lem:degenerate},
Corollary~\ref{cor:reopen}, and Proposition~\ref{prop:tunnelexclusion}. Part (III) follows from
Theorem~\ref{thm:mainA} and Corollary~\ref{cor:disconnected-even}.

\subsection{The inverting half-turn}

\begin{lemma}\label{lem:inverting}
Let $a,b$ be elliptics generating a non-elementary discrete group
$\Gamma<\PSL$. Their axes are disjoint and non-asymptotic.
If $\ell$ is their common perpendicular and $h$ is the half-turn
about $\ell$, then
\[
hah^{-1}=a^{-1},\qquad hbh^{-1}=b^{-1}.
\]
Thus $h$ normalizes $\Gamma$, and $\langle\Gamma,h\rangle$
is a discrete extension of index at most two.
Under Hypothesis~\ref{hyp:main}, with even marked orders and both
tunnels surviving in $O_2$, we have $h\in\Gamma$. Each
perpendicular foot then projects to a spherical vertex on the
image of the corresponding marked axis. A marked edge of order
at least four has a spherical endpoint in the restored graph.
\end{lemma}

\begin{proof}
Coincident, intersecting, or asymptotic marked axes would make the group
preserve a geodesic, a point of $\HH^3$, or a point at infinity,
respectively. Each alternative contradicts non-elementarity. Thus the
common perpendicular exists. Its half-turn reverses the orientation of
each marked axis and gives the two inverting identities. Consequently
it normalizes $\Gamma$, and the resulting extension, a union of at
most two translates of a closed discrete subgroup, is discrete.

Assume now the additional hypotheses of the second assertion,
and suppose $h\notin\Gamma$. Its induced isometry $\widehat h$
of $O$ is a nontrivial involution. Indeed, an isometry acting
trivially on the quotient agrees locally with one deck
transformation on the regular locus, hence agrees with it
everywhere.

At a spherical vertex, an axis of even order at least four
folds. The vertex type is $(2,2,n)$ with $n\ge4$ even, or
$(2,3,4)$, and there is only one germ of the marked order.
Such an axis projects to a single singular edge. The inverting
relations therefore make $\widehat h$ preserve each marked
edge of order at least four.

The extension $\langle\Gamma,h\rangle$ preserves the convex
core. Take a positive metric neighbourhood of that core and
truncate an invariant disjoint system of cusp neighbourhoods.
After smoothing corners equivariantly, this gives an invariant
compact core $(M'_0,P')$. The product structures on the ends
give an orientation-preserving pared homeomorphism
\[
 \Phi:(M_0,P)\longrightarrow(M'_0,P')
\]
compatible with the identifications of the fundamental groups
with $\Gamma$ and carrying each singular edge to its continuation
in $O$.

By Lemma~\ref{lem:degenerate}, the prescribed capped model has
no infinite edges or pillow vertices. Its caps are cones on
closed turnover boundaries. Cone the boundary restrictions of
$\Phi$ to obtain a weight-preserving homeomorphism from this
chosen capped model to the model $\Pi'$ obtained by capping
$M'_0$. Consequently the same marked weight change on $\Pi'$
gives an orbifold
\[
 O'_2\cong O_2\cong\Odih(r;d_+,d_-),\qquad d_+,d_-\ge2.
\]

Extend $\widehat h|_{M'_0}$ over the new caps by coning its
boundary maps. Every marked weight that changes has order at
least four, so its edge is preserved. A marked weight already
equal to two imposes no further condition. Thus the coned map
is an orientation-preserving orbifold homeomorphism of $O'_2$.
Its square is the identity because coning preserves composition,
and it is nontrivial on the interior of $M'_0$.

At least one marked order is at least four, since a discrete
group generated by two involutions is elementary. Choose that
marked edge, using $e_b$ when $p=2$. Proposition~\ref{prop:tunnelexclusion}
places its corresponding edge in $\Pi'$ on the two-bridge locus.
The coned involution stabilizes this non-tunnel edge, contradicting
Theorem~\ref{thm:AOPSYinput}(c). Hence $h\in\Gamma$.

Let $P_a\in\operatorname{Ax}(a)$ and
$P_b\in\operatorname{Ax}(b)$ be the feet of their common
perpendicular. At either foot, the rotations about the
marked axis and about $\ell$ generate a noncyclic finite
subgroup of the point stabilizer. For marked order two this
subgroup is a Klein four-group. The projected foot is
therefore an actual spherical vertex on the image of the
marked axis. For marked order at least four, it is a
spherical endpoint of the chosen marked edge.
\end{proof}

\subsection{Exact placement}

\begin{proof}[Proof of Theorem~\ref{thm:mainC}(II)(a) and of
(II)(b) when $(p,q)\ne(4,4)$]
The one-tunnel argument applies to all even $p,q$.
Normalize $1\le d_+\le d_-$. Since the tunnel weights are coprime,
the possibilities are $d_+=1<d_-$ or $2\le d_+<d_-$.

\medskip

\noindent\emph{One surviving tunnel: $d_+=1<d_-$.}
Suppose first that $K(r)=K_1\sqcup K_2$ has two components. The
surviving tunnel joins the two components and is a cut edge. Whether
its weight is odd or even, the even subgraph has first Betti number two,
with independent cycles represented by $K_1$ and $K_2$, while the tunnel
meridian has zero mod-$2$ class. Hence
\[
H_1^\orb(O_2;\mathbb F_2)\cong\mathbb F_2^2
\]
is generated by the two component-meridian classes. Since the images of
$a,b$ generate the orbifold group, their homology classes generate this
$\mathbb F_2^2$. Thus the marked edges lie on different components,
and hence are precisely the two components $K_1,K_2$.

Now suppose $K(r)$ is a knot. The tunnel feet divide it into two edges
$J_1,J_2$. If $d_-$ is odd, the tunnel cannot be marked because its
weight in $O_2$ is not $2$. Thus both marked edges lie in
$\{J_1,J_2\}$. They are either distinct, or they coincide. In the
latter case the two marked meridians are meridians of the same edge, so
their orders are equal and $p=q$. This gives the second alternative in
part~(a).

If $d_-$ is even, the even subgraph is the whole theta graph and
\[
H_1^\orb(O_2;\mathbb F_2)\cong\mathbb F_2^2.
\]
The two marked homology classes therefore generate
$\mathbb F_2^2$ and in particular are distinct. Proposition~\ref{prop:tunnelexclusion}, after changing the identification by the
symmetry of the exceptional all-$2$ theta graph if necessary, puts both
marked edges on the non-tunnel edges. Hence
\[
\{e_a,e_b\}=\{J_1,J_2\}.
\]

\medskip

\noindent\emph{Two surviving tunnels: $2\le d_+<d_-$.}
For the two-tunnel argument assume $p,q\ge4$ and
$(p,q)\ne(4,4)$. Since $d_+$ and $d_-$ are
coprime, at least one tunnel weight is odd.

Suppose first that $K(r)$ has two components. Both tunnels join the two
components. In mod-$2$ homology the odd tunnel disappears, and any even
tunnel is a cut edge of the even subgraph. Thus
\[
H_1^\orb(O_2;\mathbb F_2)\cong\mathbb F_2^2
\]
is freely generated by the two component classes. Consequently the
marked edges lie on different components.

At either endpoint of a marked edge, the other knot arc on the same
component is unmarked and has weight $2$. Hence the two endpoint types
of $e_a$ are $(p,2,d_+)$ or $(p,2,d_-)$, and similarly for $e_b$.
By Lemma~\ref{lem:inverting}, one endpoint of $e_a$ is spherical, so
for some $d\in\{d_+,d_-\}$,
\[
\frac1p+\frac12+\frac1d>1,
\qquad\text{equivalently}\qquad
\frac1p+\frac1d>\frac12.
\]
Since $p\ge4$, this forces $d\le3$. If $d_+\ne2$, then
$d_+=3$ and $d_->3$. The only possible spherical foot of $e_a$ is
then the $d_+$-foot, and the displayed inequality forces $p=4$; the
same argument for $e_b$ forces $q=4$, contrary to
$(p,q)\ne(4,4)$. Therefore $d_+=2$, and coprimality makes $d_-$ odd.
The standard horizontal involution puts the marks on $J_{1,1}$ and
$J_{i,2}$ as stated. For $i=1,2$ the two placements correspond,
respectively, to slopes $r$ and $r+1$ in the convention of
\cite[Section~9, Subcase~3.1]{AOPSY}.

It remains to suppose that $K(r)$ is a knot. Write its four arcs
cyclically as in Theorem~\ref{thm:mainC}, and use the model symmetry
from Theorem~\ref{thm:AOPSYinput}(c) to normalize
$e_a=e_1$. Proposition~\ref{prop:tunnelexclusion} puts $e_b$ on
$K(r)$ as well.

Assume first that exactly one tunnel has even weight. The even
subgraph is $K(r)$ together with that tunnel, a theta graph, and
\[
H_1^\orb(O_2;\mathbb F_2)\cong\mathbb F_2^2.
\]
Hence $e_b\ne e_a$. If $e_b$ is adjacent to $e_a$ across the odd
tunnel, Lemma~\ref{lem:mod2weighted} identifies their homology classes,
which is impossible. Suppose instead that they are adjacent across the
even tunnel. Their common endpoint has type $(p,q,d_{\rm ev})$, which
is not spherical because
\[
\frac1p+\frac1q+\frac1{d_{\rm ev}}
\le\frac14+\frac14+\frac12=1.
\]
Lemma~\ref{lem:inverting} therefore forces the other endpoint of each
marked edge to be spherical. Those endpoint types are
\[
(p,2,d_{\rm odd}),\qquad(q,2,d_{\rm odd}),
\]
so
\[
\frac1p+\frac1{d_{\rm odd}}>\frac12,\qquad
\frac1q+\frac1{d_{\rm odd}}>\frac12.
\]
Since $d_{\rm odd}\ge3$ is odd, these inequalities force
$d_{\rm odd}=3$ and $p=q=4$, contrary to the hypothesis. Hence $e_b$
is opposite $e_a$, namely
\[
e_b=e_3.
\]

The endpoints of $e_a$ are now $(p,2,d_+)$ and $(p,2,d_-)$, so the
sphericity inequality holds for one of the two tunnel weights. The inequality forces $d_+\le3$. If
$d_+=3$, the preceding argument applied to both marked edges yields $p=q=4$. Therefore
$d_+=2$, and coprimality makes $d_-$ odd.

Finally suppose both tunnel weights are odd. A marked edge lies on
$K(r)$. At a spherical endpoint of $e_a$, the three incident weights
have the form $(p,c,d)$, where $c$ is either $2$ or $q$, and $d\ge3$.
Thus
\[
\frac1p+\frac1c+\frac1d
\le
\frac1p+\frac12+\frac13.
\]
If $p\ge6$, the right-hand side is at most $1$, contradicting
sphericity. Hence $p=4$. The same argument gives $q=4$, again excluded.
Thus the both-odd case cannot occur.
\end{proof}

\section{The exceptional orders \texorpdfstring{$(4,4)$}{(4,4)}}\label{sec:exceptional}

Throughout this section assume the hypotheses of Theorem~\ref{thm:mainC},
with both tunnels surviving and $p=q=4$. We retain the
inverting half-turn $h$ from Lemma~\ref{lem:inverting}. Its axis is the
common perpendicular $\ell$, and $h\in\Gamma$. Thus every point of
$\ell$ projects into the finite singular locus. In the interior of an edge
its label must be even, since its local cyclic stabilizer contains $h$.
Write $P_a,P_b$ for the perpendicular feet. These are distinct points
of $\HH^3$, although their images in the quotient need not be distinct.

\subsection{Continuation of a fixed axis}

At a spherical vertex, incident edge germs are orbits of rays under
the local rotation group. The two directions of a lifted axis may
therefore project to the same germ. Continuing the parametrized
geodesic then folds its projected path; it does not end the geodesic.

\begin{lemma}\label{lem:44tracklocal}
Let the axis of an involution pass through a spherical vertex of an
orientable hyperbolic $3$-orbifold.
\begin{enumerate}[(i)]
\item At a vertex of type $(2,3,4)$, an involution axis perpendicular
to an order-$4$ axis can be either an order-$2$ axis or an order-$4$
axis. Its two rays project to the same germ, of label $2$ or $4$
respectively.
\item At a vertex of type $(2,2,m)$, an involution axis along a
side order-$2$ germ continues to the other order-$2$ germ if $m$ is
odd. If $m$ is even, both rays project to the same germ. For $m=2$
each of the three germs has this latter property.
\end{enumerate}
In the cases where the rays project to the same germ, the projected
axis folds at the vertex. It does not terminate there.
\end{lemma}

\begin{proof}
For (i), realize the octahedral group as the determinant-one signed
permutation matrices on $\mathbb R^3$, and take the order-$4$ axis to
be the third coordinate axis. The perpendicular involution axes are
the first two coordinate axes and the two diagonal axes in their plane.
The former have generic stabilizer of order $4$; the latter have
generic stabilizer of order $2$. For each such axis, its two directions
are in the same orbit under the octahedral group.

For (ii), realize $D_m$ as rotations about the third coordinate axis
and half-turns about the axes at angles $j\pi/m$ in the horizontal
plane. The opposite directions on such an axis belong to distinct
direction orbits precisely when $m$ is odd: the rotation generator
changes $j$ by $2$, whereas reversal changes it by $m$. This gives
the two side germs for odd $m$, and the stated folding for even $m$.
The Klein-four case is immediate in coordinate axes.
\end{proof}

The continuation rule is applied to the parametrized geodesic $\ell$,
not just to its image as a subset. Along an open edge the projection is
a local isometry. A compact subsegment cannot reach a deleted Euclidean
or hyperbolic vertex; a change of direction on the graph can occur only
through the local quotient at a spherical vertex.

\subsection{The additional configurations}

If a marked $4$-edge has an unmarked $2$-edge and a tunnel of weight
$t$ at a spherical end, then
\[
 \frac14+\frac12+\frac1t>1,
 \qquad\text{hence }t=2\text{ or }3.
\]
A vertex with both marked $4$-germs is never spherical. Consequently,
apart from the placement of Theorem~\ref{thm:mainC}(II)(b), the following
are the possible knot configurations. We use the cyclic edge convention
of that theorem; in this table the tunnels are not necessarily ordered
by weight.
\[
\begin{array}{c|cccc|cc|l}
 & e_1&e_2&e_3&e_4&\tau_+&\tau_-&\text{restrictions}\\\hline
 A_m&4&4&2&2&m&3&m\ge2\text{ even},\ (m,3)=1\\
 B&4&2&4&2&3&d&d\ge4,\ (d,3)=1\\
 C&4&2&2&2&3&d&d\ge5\text{ odd},\ (d,3)=1\\
 D&4&2&2&4&3&d&d\ge5\text{ odd},\ (d,3)=1.
\end{array}
\]
In $C$ both marks are on $e_1$, and in the other rows they are
on the two $4$-edges. Each row has the embedding
$K(r)\cup\tau_+\cup\tau_-$. Interchanging the tunnels identifies
the incidence patterns of $A_m$ and $D$. The distinction is the
parity of the tunnel opposite the weight-$3$ tunnel. At a
$(2,2,m)$ vertex the inverting axis folds when $m$ is even and
continues through when $m$ is odd. This is why the axis-tracing
argument below excludes $A_m$, while $D$ requires the distance
estimate.

\begin{lemma}\label{lem:44list}
The table is exhaustive for additional knot configurations. For a
two-component link, the only additional possibility has tunnel weights
$3,d$ with $d\ge4$ coprime to $3$, and one marked $4$-arc on each
component.
\end{lemma}

\begin{proof}
For a link with two components, mod-$2$ homology puts the marks on
different components. Each marked arc has an unmarked $2$-arc at
each end. Its spherical-end condition forces the smaller tunnel weight
to be $2$ or $3$. Weight $2$ gives the usual alternative, and weight
$3$ gives the stated additional case.

For a knot, if exactly one tunnel is even, the two marks are distinct
and cannot be adjacent across the odd tunnel, by mod-$2$ homology.
If they are opposite, their spherical ends force a tunnel of weight
$2$ or $3$. Weight $2$ gives the usual alternative; otherwise the
configuration is $B$. If they are adjacent across the even tunnel,
their common end is nonspherical and their other ends force the odd
tunnel weight to be $3$. This is $A_m$ for an arbitrary even $m$
coprime to $3$.

If both tunnels are odd, the smaller must be $3$. The marks may
coincide, giving $C$, or be opposite, giving $B$. For adjacent marks,
their common end is nonspherical. Their other ends must therefore be
at the weight-$3$ tunnel; equivalently, they are adjacent across the
other tunnel. This gives $D$.
\end{proof}

\begin{proposition}\label{prop:44AB}
The configurations $A_m$, $B$, and $C$ do not occur. The additional
two-component configuration in Lemma~\ref{lem:44list} does not occur.
\end{proposition}

\begin{proof}
Put $v_{ij}=e_i\cap e_j$ in the knot case. All assertions concern the
projection of the segment $[P_a,P_b]\subset\ell$, with continuation
governed by Lemma~\ref{lem:44tracklocal}.

In $A_m$, the vertex $v_{12}$ has type $(4,4,m)$ and is nonspherical;
$v_{23},v_{41}$ have type $(4,2,3)$; and $v_{34}$ has type $(2,2,m)$.
The two feet must project to the two different octahedral vertices.
At either such vertex, the fixed axis can leave on the marked $4$-edge
or on the adjacent $2$-edge. The former leads to the nonspherical
vertex $v_{12}$. The latter leads to $v_{34}$, where it folds back
along the same edge because $m$ is even. At the octahedral vertex it
again folds along the same $2$-germ. Thus this continued axis never
reaches the other octahedral vertex. Neither possibility connects the
two feet.

In $B$, the only spherical ends of the marked edges are $v_{12}$ and
$v_{34}$, both of type $(4,2,3)$. The other two vertices have type
$(4,2,d)$ and are nonspherical for $d\ge4$. From either octahedral
foot, both possible outgoing even germs lead directly to a nonspherical
end. No segment of the fixed axis can connect the two feet. In the
two-component additional configuration, the same argument applies on
each component: both even germs at its octahedral vertex lead to its
nonspherical $(4,2,d)$ vertex, and the weight-$3$ tunnel cannot carry
the axis of $h$.

In $C$, the only spherical end of the marked edge is $v_{12}$ of
type $(4,2,3)$. Both feet must project there. The other vertex types
are
\[
 v_{23}:(2,2,d),\qquad v_{34}:(2,2,3),\qquad
 v_{41}:(4,2,d).
\]
The last is nonspherical. An axis leaving $v_{12}$ on the $4$-germ
goes directly towards that nonspherical end. An axis leaving on the
$2$-germ follows $e_2,e_3,e_4$ in order: at $v_{23}$ and $v_{34}$
the odd dihedral rule continues it to the other $2$-germ. It again
runs towards $v_{41}$ without a further fold. On either side of the
initial octahedral point, the lifted geodesic therefore never returns
to a point projecting to $v_{12}$. Hence it contains at most one such
point. This contradicts the distinctness of $P_a$ and $P_b$.
\end{proof}

\Needspace{5\baselineskip}
\subsection{A uniform bound for a generating pair}

\begin{lemma}[Order-four trace rectangle]\label{lem:44-rectangle}
Let $A,B\in\mathrm{SL}_2(\mathbb C)$ have trace $\sqrt2$, and let
$a,b$ be their projective images. Suppose that $\langle a,b\rangle$
is discrete and the marked epimorphism
$C_4*C_4\to\langle a,b\rangle$ has nontrivial kernel. Then
\[
 -2<\Re\tr(AB)<4,\qquad |\Im\tr(AB)|<\sqrt3.
\]
\end{lemma}

\begin{proof}
Put $z=\tr(AB)$ and $\gamma=\tr[A,B]-2$. The Fricke identity
\[
 \tr[A,B]=x^2+y^2+z^2-xyz-2,
 \qquad x=\tr A,\quad y=\tr B,
\]
gives $\gamma=z(z-2)$. For $p=q=4$, the parameter equation of
\cite[Theorem~1.1]{EGMS} is $\rho(\rho-2)=\gamma$, whose roots
are $z$ and $2-z$. Let $\Omega_{4,4}$ be the open polygon defined in
\cite[Definition~2.6]{EGMS}. Theorem~1.1 of that paper makes
the marked representation discrete and faithful whenever
$\rho\notin\Omega_{4,4}$. By parts (1)--(2) of that definition,
the polygon is contained in
\[
 \{\rho\in\mathbb C:-2<\Re\rho<4,\quad |\Im\rho|<\sqrt3\}.
\]
Indeed, its vertical supporting lines are $\Re\rho=4,-2$, and
its horizontal supporting lines are $\Im\rho=\pm\sqrt3$.
The contrapositive proves the bounds. Replacing $b$ by $b^{-1}$
sends $z$ to $\tr(AB^{-1})=2-z$, and preserves this rectangle.
\end{proof}

\begin{proposition}\label{prop:44bound}
Let $a,b$ be primitive order-$4$ elliptics which generate a
non-elementary discrete group not freely generated as $C_4*C_4$.
Choose lifts $A,B$ of trace $\sqrt2$, and let $\ell_0$ be the
distance between their axes. Then
\[
 |\tr[A,B]-2|\le\sqrt{133}.
\]
In particular,
\begin{equation}\label{eq:44-upper}
 \ell_0\le L_*:=\operatorname{arsinh}(133^{1/4}),
 \qquad \cosh\ell_0\le\sqrt{1+\sqrt{133}}.
\end{equation}
\end{proposition}

\begin{proof}
Put $z=\tr(AB)$ and $\gamma=\tr[A,B]-2=z(z-2)$.
By Lemma~\ref{lem:44-rectangle},
$t=(\Re z-1)^2\le9$ and $u=(\Im z)^2\le3$. Hence
\[
 |\gamma|^2=(t+u+1)^2-4t
 \le(t+4)^2-4t=t^2+4t+16\le133.
\]
For oriented axes, the complex distance
$\delta=\ell_0+i\theta$ records both their distance and the angle
between their directions after parallel transport along the common
perpendicular. The trace--distance identity for primitive
order-$4$ rotations is $\gamma=\sinh^2\delta$. Consequently
\[
 \sinh^2\ell_0\le
 \sinh^2\ell_0+\sin^2\theta=|\gamma|\le\sqrt{133}.
\]
This proves both bounds in \eqref{eq:44-upper}.
\end{proof}

\subsection{The exchanging half-turn and the midpoint stabilizer}

The following construction applies to any equal primitive orders.

\begin{lemma}[Exchanging half-turn]\label{lem:44exchange}
Let $a,b$ be primitive elliptic meridians of the same order
$n\ge3$ generating a non-elementary discrete group $\Gamma<\PSL$.
Write $[P_a,P_b]$ for the common perpendicular segment of their
axes, and let $h$ be the inverting half-turn of
Lemma~\ref{lem:inverting}. There is a half-turn $f$ such that
\[
 faf^{-1}=b,\qquad fbf^{-1}=a,\qquad fh=hf.
\]
It fixes the midpoint $M$ of the common perpendicular segment
$[P_a,P_b]$, and its axis is perpendicular to $\ell=\ax h$ there.
The extension $\widehat\Gamma=\langle\Gamma,f\rangle$ is discrete
of index at most two over $\Gamma$. If $n\ge4$ and $a,b$ are
meridians of distinct singular edges, then $f\notin\Gamma$.
\end{lemma}

\begin{proof}
Choose lifts $A,B$ with $\tr A=\tr B=x=2\cos(\pi/n)$, and put
\[
 F=A+B-xI.
\]
Cayley--Hamilton gives
\[
 \tr F=0,\qquad FA=BF,\qquad FB=AF,
 \qquad F^2=(\tr(AB)-2)I.
\]
The matrix $F$ is invertible. Indeed, $F=0$ would give
$B=A^{-1}$. Otherwise, a nonzero kernel of $F$ would be a common
invariant line for $A$ and $B$, by the two intertwining identities.
Both alternatives contradict non-elementarity. After multiplication
by a scalar, $F$ represents an element $f$ of $\PSL$.
Its trace is zero, so it is a half-turn, and the intertwining
identities give the exchanging relations.

For a lift $H$ of $h$, the nonzero common trace $x$ fixes the signs
in the lifted inverting relations. Thus
$HAH^{-1}=A^{-1}=xI-A$ and $HBH^{-1}=B^{-1}=xI-B$.
Consequently $HFH^{-1}=-F$, so $f$ and $h$ commute in $\PSL$.
The isometry $f$ interchanges the marked axes and hence the feet
of their unique common perpendicular. It therefore fixes their
midpoint and reverses $\ell$. A half-turn reversing a geodesic
has axis perpendicular to it at the fixed point.

The exchanging relations show that $f$ normalizes $\Gamma$.
The extension is a union of at most two translates of a discrete
subgroup, hence is discrete.
Finally suppose $n\ge4$. At a spherical vertex an order-$n$
axis has a single germ of that order. The possible triples are
$(2,2,n)$ and, for $n=4,5$, $(2,3,n)$.
Its projected axis consequently lies on a single singular edge.
Conjugate elements of $\Gamma$ have axes with the same image
in $O$, so distinct marked edges imply $f\notin\Gamma$.
\end{proof}

\begin{lemma}[Localization of the midpoint]\label{lem:Dmidpoint}
Suppose that configuration $D$ occurs, with its marked generating
pair on $e_1,e_4$. For the midpoint $M$ in
Lemma~\ref{lem:44exchange}, one has
\[
 \Gamma_M\cong D_d,\qquad
 \widehat\Gamma_M\cong D_{2d}.
\]
The cyclic subgroup of order $2d$ in $\widehat\Gamma_M$ has an
axis perpendicular to $\ell$ at $M$.
\end{lemma}

\begin{proof}
Let $\pi:\HH^3\to O$ be the quotient map, and denote the spherical
vertices of $D$ by
\[
 v_{12}:(4,2,3),\qquad
 v_{23}:(2,2,d),\qquad
 v_{34}:(2,4,3).
\]
The remaining formal vertex $v_{41}$ has type $(4,4,d)$ and is
hyperbolic; it is not a point of $O$.

Since the marked edges are distinct, $f\notin\Gamma$. Its induced
involution $\bar f$ on $O$ exchanges $e_1,e_4$. Their unique
spherical endpoints are therefore exchanged:
$v_{12}\leftrightarrow v_{34}$. Incidence and edge orders then force
$e_2\leftrightarrow e_3$ and preserve each tunnel. In particular,
no interior point of any of the four even edges is fixed by $\bar f$,
and neither octahedral vertex is fixed.

On the other hand, $fM=M$, so $\bar f$ fixes $\pi(M)$.
Since $h\in\Gamma$ and $hM=M$, the stabilizer of $\pi(M)$
contains an involution. Thus $\pi(M)$ is neither a regular point
nor an interior point of either odd-labelled tunnel. The only
remaining actual singular vertex is $v_{23}$. Hence
$\pi(M)=v_{23}$ and $\Gamma_M\cong D_d$.

Normality of $\Gamma$ in $\widehat\Gamma$, together with
$f\in\widehat\Gamma_M\setminus\Gamma_M$, gives
\[
 1\longrightarrow\Gamma_M\longrightarrow\widehat\Gamma_M
 \longrightarrow C_2\longrightarrow1.
\]
The latter stabilizer is finite, of order $4d$. Since $d$ is odd,
the rotation subgroup $C_d$ of $D_d$ is characteristic: it consists
of all elements of odd order. Consequently
$\widehat\Gamma_M$ preserves its rotation axis. A finite subgroup
of $\mathrm{SO}(3)$ preserving an unoriented axis is cyclic or
dihedral. It contains $D_d$, so it is dihedral, and its order makes
it $D_{2d}$. Finally, $h\in D_d$ is an involution and $d$ is odd;
thus $h$ is a side half-turn, whose axis is perpendicular to the
rotation axis at $M$.
\end{proof}

\subsection{Exclusion of Family D}

\begin{theorem}\label{thm:Dexcluded}
Configuration $D$ is impossible for every odd $d\ge5$ coprime to
$3$, independently of the two-bridge slope.
\end{theorem}

\begin{proof}
Suppose otherwise. Put
$\ell_0=\operatorname{dist}(\ax a,\ax b)>0$ and retain the
midpoint $M$. By Lemma~\ref{lem:Dmidpoint},
$\widehat\Gamma$ contains a primitive rotation $c$ of order $2d$
whose axis is perpendicular to $\ell$ at $M$. The axis of $a$ is
perpendicular to $\ell$ at $P_a$. Their common perpendicular is
therefore $[P_a,M]$, of length $\ell_0/2$.

The subgroup $\langle c,a\rangle$ is discrete and non-elementary.
Its two elliptics have orders greater than two and have disjoint,
non-asymptotic axes. A finite elementary group would fix a point
of $\HH^3$, forcing these axes to meet. An infinite elementary
group fixing an ideal point would force a common ideal endpoint.
If it preserved an unordered pair of ideal points, both elliptics
would fix that pair pointwise, since their orders exceed two.
Their axes would then coincide. Each case contradicts the geometry
of the two axes.

Take lifts $C,A$ to $\mathrm{SL}_2(\mathbb C)$ and put
$s=\sin(\pi/(2d))$. J\o rgensen's inequality
\cite[Lemma 1]{Jorgensen76} gives
\[
 1\le |\tr^2C-4|+|\tr[C,A]-2|.
\]
If $\ell_0/2+i\theta$ is the complex distance between the axes,
the trace--distance formula gives
\[
 |\tr^2 C-4|=4s^2,\qquad
 |\tr[C,A]-2|
 =2s^2\bigl(\sinh^2(\ell_0/2)+\sin^2\theta\bigr)
 \le2s^2\cosh^2(\ell_0/2).
\]
Consequently
\[
 \cosh^2(\ell_0/2)\ge\frac1{2s^2}-2,
 \qquad
 \cosh\ell_0\ge\frac1{s^2}-5.
\]
Since $d\ge5$, one has
$s^2\le\sin^2(\pi/10)=(3-\sqrt5)/8$, and hence
\begin{equation}\label{eq:D-lower}
 \cosh\ell_0\ge 1+2\sqrt5.
\end{equation}
But Proposition~\ref{prop:44bound} gives
$\cosh\ell_0\le\sqrt{1+\sqrt{133}}$.
These bounds are incompatible: the latter number is less than $4$,
whereas $1+2\sqrt5>5$.
\end{proof}

\begin{corollary}\label{thm:44reduction}
Under the hypotheses of Theorem~\ref{thm:mainC}, suppose $p=q=4$ and
both tunnels survive. Then their weights are $2,m$ with $m\ge3$
odd. For a two-component link the marks lie on different components;
for a knot they are opposite arcs. Thus the placement conclusion of
Theorem~\ref{thm:mainC} holds for all even $p,q\ge4$.
\end{corollary}

\begin{proof}
Lemma~\ref{lem:44list} lists the additional configurations.
Proposition~\ref{prop:44AB} excludes $A_m,B,C$ and the additional
two-component configuration. Theorem~\ref{thm:Dexcluded} excludes
$D$. The one-tunnel argument was independent of the exception.
\end{proof}

Corollary~\ref{thm:44reduction} completes the proof of
Theorem~\ref{thm:mainC}(II).

\section{An order-two generator}\label{sec:ordertwo}

The proofs of Theorem~\ref{thm:mainC}(I) and (II)(a) apply to all
even marked orders. We first record an index-two construction
for an involution and an elliptic of arbitrary order $q\ge3$.

\begin{lemma}\label{lem:indextwo}
Assume Hypothesis~\ref{hyp:main}, with $p=2$ and $q\ge3$.
Put $b'=aba$ and $\Gamma_0=\langle b,b'\rangle$. Then $\Gamma_0$ is
the normal closure of $b$ in $\Gamma$ and has index at most two.
It is discrete and non-elementary, $b,b'$ are primitive meridians
of order $q$ in $\Gamma_0$, and the marked epimorphism
$C_q*C_q\twoheadrightarrow\Gamma_0$ is not injective.
\end{lemma}

\begin{proof}
Conjugation by $a$ exchanges $b$ and $b'$, so $\Gamma_0$ is normal
and is the normal closure of $b$. Its quotient in $\Gamma$ is
generated by the image of $a$, hence has order at most two.
Discreteness and non-elementarity pass to a subgroup of finite index.
The subgroup of $\Gamma$ fixing $\ax b$ pointwise is
$\langle b\rangle$. Its intersection with $\Gamma_0$ is unchanged,
which proves primitivity in the cover; conjugation by $a$ gives the
same assertion for $b'$.

Let $F=C_2*C_q=\langle A,B\mid A^2=B^q=1\rangle$ and let
$F_0$ be the kernel of $F\to C_2$ sending $A$ to the nonidentity
element and $B$ to the identity. The subgroup theorem for free
products gives
\[
 F_0=\langle B,ABA\rangle\cong C_q*C_q.
\]
Write $N$ for the kernel of the marked map $F\to\Gamma$.
If $F_0\to\Gamma_0$ were injective, then $N\cap F_0=1$ and $N$
would inject into $F/F_0$, so $N$ would be finite. A free product
of two nontrivial finite groups has no nontrivial finite normal
subgroup: such a subgroup fixes a vertex of the Bass--Serre tree,
its fixed subtree is invariant by normality, and minimality and
trivial edge stabilizers force the subgroup to be trivial
\cite{Serre}. Thus $N=1$, contrary to Hypothesis~\ref{hyp:main}.
\end{proof}

For the remainder of this section assume $p=2$ and $q\ge4$ is even.
The case $q=2$ is excluded because a discrete group generated
by two involutions is elementary.

We call the union of the closed edges traversed by the image of
$\ax a$ the \emph{chain} of $a$. Every edge of this chain has
weight two. Indeed, a rotation fixing an open segment of $\ax a$
fixes the whole axis, whose pointwise stabilizer is $\langle a\rangle$.
At a spherical vertex of type $(2,2,n)$ with $n$ odd, the chain
continues between the two order-two germs. At every other spherical
vertex that occurs here it folds, by Lemma~\ref{lem:44tracklocal}.
Each of its edges represents $a$ as a meridian with a suitable base
path. Consequently Proposition~\ref{prop:tunnelexclusion} applies
to each edge of the chain: when both tunnels survive, the chain
lies on $K(r)$. 
\begin{proposition}\label{prop:twofour}
Assume $p=2$, $q\ge4$ is even, and both tunnels survive. Then,
after normalization,
\[
 d_+=2,\qquad d_-=m\ge3\quad\hbox{with $m$ odd}.
\]
If $K(r)$ has two components, the chain of $a$ is the entire
component not containing $e_b$. If $K(r)$ is a knot, use the cyclic
labelling of Theorem~\ref{thm:mainC}(II)(b) with $e_b=e_1$.
Then the chain of $a$ is $e_2\cup e_3$. In particular, its marked
edge may be chosen to be the opposite edge $e_3$.
\end{proposition}

\begin{proof}
Let $h$ be the inverting half-turn of Lemma~\ref{lem:inverting}.
Then $h\in\Gamma$, the foot $P_b$ projects to a spherical endpoint
of $e_b$, and the foot $P_a$ projects to a spherical vertex on the
chain of $a$ whose stabilizer contains $\langle a,h\rangle\cong D_2$.
All link edges other than $e_b$ have weight two. Thus a spherical
endpoint of $e_b$ has type $(q,2,t)$, where $t$ is a tunnel weight,
and
\[
 \frac1q+\frac12+\frac1t>1.
\]
This forces $t=2$ for $q\ge6$, and $t\in\{2,3\}$ for $q=4$.
Since the two tunnel weights are coprime and exceed one, either
$d_+=2$ and $d_-$ is odd, or $q=4$, $d_+=3$, and
$d_-=d\ge4$ with $(3,d)=1$.

\medskip
\noindent\emph{Two components.}
The mod-$2$ homology computation in Section~\ref{sec:connected}
gives two independent component-meridian classes. Hence the chain
of $a$ lies on the component $K_1$ not containing $e_b$.
If $d_+=2$ and $d_-=m$ is odd, the chain continues through the
$(2,2,m)$ vertex and folds at the $(2,2,2)$ vertex, and therefore
traverses both arcs of $K_1$.

Suppose $d_+=3$. The vertices on $K_1$ have types $(2,2,3)$ and
$(2,2,d)$. Only the latter can contain $D_2$, and only when $d$
is even. Thus $P_a$ projects to the $(2,2,d)$ vertex $w$.
The axis of $h$ at $w$ is either the rotation axis or a side axis
of $D_d$. In the first case its image follows $\tau_-$ towards
the nonspherical $(4,2,d)$ vertex on the other component; a finite
subsegment cannot reach that deleted vertex. In the second case
its image follows $K_1$, continuing at the odd $(2,2,3)$ vertex
and folding whenever it returns to $w$. It therefore stays on
$K_1$. Neither case reaches the octahedral vertex on the other
component to which $P_b$ projects. This contradicts the fact that
$P_a$ and $P_b$ lie on the axis of $h$.

\medskip
\noindent\emph{A knot.}
Use the standard symmetry to put $e_b=e_1$ and write
$v_{ij}=e_i\cap e_j$ for consecutive indices. If $d_+=2$ and
$d_-=m$ is odd, Lemma~\ref{lem:mod2weighted} gives
\[
 [e_1]=[e_4],\qquad [e_2]=[e_3],
\]
with the two displayed classes independent. Since $a,b$ generate,
the chain of $a$ contains $e_2$ or $e_3$. It continues through
$v_{23}$, of type $(2,2,m)$, and folds at $v_{12}$ and $v_{34}$,
of types $(q,2,2)$ and $(2,2,2)$. Its image is therefore precisely
$e_2\cup e_3$.

It remains to exclude $q=4$ and tunnel weights $3,d$, where
$d\ge4$ and $(3,d)=1$. The successive vertex types are
\[
 (4,2,3),\qquad(2,2,d),\qquad(2,2,3),\qquad(2,4,d).
\]

Apply Lemma~\ref{lem:indextwo}. If $d$ is odd, both tunnels
have odd weight, so the even subgraph is the single cycle
$K(r)$. Lemma~\ref{lem:mod2weighted} gives
\[
H_1(\Gamma;\mathbb F_2)\cong\mathbb F_2,
\]
generated by the class of $b$, a meridian of $e_1$.
Hence every homomorphism $\Gamma\to C_2$ killing $b$ is
trivial. Since $\Gamma_0$ is the normal closure of $b$ and
has index at most two, $\Gamma_0=\Gamma$. Thus $\Gamma$
is generated by the order-four meridians $b,aba$ of the
same edge $e_1$. The original surgery changes only the
weight of $e_1$, since $a$ already has order two. Surgery
for the pair $(b,aba)$ lowers this same weight from four
to two, and hence gives the same chosen model $O_2$.
By Lemma~\ref{lem:indextwo} the
$(4,4)$ placement theorem applies, and it forbids coincident
marked edges when both tunnels survive.

Suppose $d$ is even. In this case the mod-$2$ relations are
\[
 [e_1]=[e_2],\qquad [e_3]=[e_4],\qquad
 [\tau_-]=[e_1]+[e_3],\qquad [\tau_+]=0,
\]
and the two link-edge classes are independent. There is therefore
a nontrivial homomorphism to $C_2$ killing $b$. By
Lemma~\ref{lem:indextwo}, $[\Gamma:\Gamma_0]=2$. The edges whose
meridians have nontrivial image are exactly
\[
 C=e_3\cup e_4\cup\tau_-.
\]
Consider the corresponding orbifold cover and extend it over the
capping balls. On underlying spaces it is branched along $C$.
An unbranched edge has two lifts of the same weight; a branched
edge has one lift with half its weight, deleted when that weight
is one. The vertices lift as follows:
\begin{center}\small
\begin{tabular}{@{}llll@{}}
\toprule
Vertex & Type & Branch germs & Lift\\
\midrule
$v_{12}$ & $(4,2,3)$ & none & two $(4,2,3)$ vertices\\
$v_{23}$ & $(2,2,d)$ & $e_3,\tau_-$ & one $(2,2,d/2)$ vertex\\
$v_{34}$ & $(2,2,3)$ & $e_3,e_4$ & a smooth order-three edge\\
$v_{41}$ & formal $(2,4,d)$ & $e_4,\tau_-$ & one formal $(4,4,d/2)$ vertex\\
\bottomrule
\end{tabular}
\end{center}
The spherical rows follow by taking the kernels of the local
homomorphisms to $C_2$: at $v_{23}$ the kernel is $D_{d/2}$,
whereas at $v_{34}$ it is $C_3$. For the last row, take the double cover of
the boundary turnover: the unbranched order-four point has two
preimages of order four, the branched order-two point becomes
regular, and the remaining cone order is $d/2$. The underlying
double cover of $S^2$ is branched at two points, so
Riemann--Hurwitz gives $\chi(\widetilde S)=2\chi(S^2)-2=2$.
Thus its underlying surface is a sphere. The reciprocal sum of
the remaining cone orders is $1/2+2/d$, so the turnover is
Euclidean for $d=4$ and hyperbolic for $d>4$.

This construction is a fully capped model of
$\mathbb H^3/\Gamma_0$, with its parabolic data pulled back
from the original core. Lemma~\ref{lem:background} identifies
its closed ambient manifold with $S^3$.

Write $v_+,v_-$ for the two lifts of $v_{12}$, write $u$ for
the lift of $v_{23}$, and write $w$ for the formal vertex above
$v_{41}$. The lifted capped graph is the complete graph on
$v_+,v_-,u,w$. Its edges $uv_\pm$ have weight two, its edges
$wv_\pm$ have weight four, and $uw$ has weight $d/2$.
The two lifts of $\tau_+$ concatenate at the smooth lift of
$v_{34}$ to give the edge $v_+v_-$ of weight three.
The marked edges are $wv_+,wv_-$, the two lifts of $e_1$.
The deck transformation represented by $a$ exchanges them,
and their meridians are $b,aba$.

Lemma~\ref{lem:indextwo} makes the $(4,4)$ classification
applicable to this chosen model. A one-tunnel model has two
vertices, whereas the present graph has four. In a two-tunnel
model the marked edges are disjoint, by
Corollary~\ref{thm:44reduction}, whereas $wv_+$ and $wv_-$
share the vertex $w$. This excludes $d_+=3$ and completes
the proof.
\end{proof}

\section{Geometric uniqueness of the spatial models}\label{sec:rigidity}

\begin{proposition}\label{prop:pared-rigidity}
Let $\Gamma_1,\Gamma_2<\PSL$ be non-elementary geometrically finite
groups. Suppose that an orientation-preserving orbifold homeomorphism
\[
 \Phi:(M_0^1,P^1)\longrightarrow(M_0^2,P^2)
\]
identifies their pared compact cores and their cusp data. If every
component of $\Omega(\Gamma_1)/\Gamma_1$ is a hyperbolic turnover,
with punctures allowed, then $\Phi_*$ is induced by conjugation with
a M\"obius transformation, up to the inner automorphism determined
by the choice of base paths. An orientation-reversing identification
is induced instead by an anti-M\"obius transformation.
\end{proposition}
\begin{proof}
Put $O_i=\HH^3/\Gamma_i$ and $(M,P)=(M_0^1,P^1)$.
Use the identity and $\Phi$ to regard $O_1$ and $O_2$ as marked
geometrically finite structures on the same oriented pared orbifold
$(M,P)$, with its prescribed cusp data.
In the deformation space of \cite[Definition~8.4.2]{BP},
structures are marked by proper homotopy classes of
orientation-preserving homeomorphisms of $(M,P)$. Two structures
are equivalent when an isometry realizes their marking change.
For fixed parabolic locus $P$, this space is parametrized by the
Teichm\"uller space of the nonparabolic boundary, with its boundary
circles treated as punctures
\cite[Proposition~8.4.4]{BP}.
Each turnover has a unique marked conformal structure, so the
deformation space is a singleton. The empty-boundary case is
\cite[Remark~8.4.5]{BP}. 
Hence there is an orientation-preserving
orbifold isometry $I:O_1\to O_2$ representing the proper homotopy
class determined by $\Phi$. Lifting $I$ to $\HH^3$ gives a
M\"obius transformation whose conjugation induces $\Phi_*$,
up to the inner automorphism arising from base paths.
If $\Phi$ reverses orientation, first compose it on the target with
the isometry induced by complex conjugation, replacing $\Gamma_2$
by its complex conjugate. The preceding argument applies to this
orientation-preserving marking. Composing the resulting lift with
complex conjugation gives the required anti-M\"obius transformation.
\end{proof}

\begin{corollary}\label{cor:model-holonomy}
Under the hypotheses of Theorem~\ref{thm:mainC}, let $Q$ be the
pared orbifold obtained by restoring the marked weights and reopening
the prescribed vertices of its classified spatial model. Then
$\Gamma\cong\pi_1^{\orb}(Q)$. Every orientation-compatible
complete hyperbolic realization of $Q$, with the prescribed
cusp data, has holonomy M\"obius conjugate to $\Gamma$. The
identification preserves the conjugacy classes of the marked edge
meridians, up to inversion.
\end{corollary}

\begin{proof}
The spatial identification and the prescribed realization give an
orbifold homeomorphism of pared cores, hence the group isomorphism
and the statement about meridians. The marked meridians in a second
complete realization have the same primitive local orders and the
same nontrivial marked relation. Its faithful holonomy is
non-elementary, since the abstract group is not virtually abelian.
Proposition~\ref{thm:rigid} makes both realizations geometrically
finite, and Lemma~\ref{lem:rigid-boundary-turnover} gives turnover
conformal boundary. Proposition~\ref{prop:pared-rigidity} now applies.
\end{proof}

\begin{proof}[Proof of Theorem~\ref{thm:classification}]
The subgroup of $\Gamma$ fixing $\ax g_1$ pointwise is finite
cyclic. Choose its once-around rotation generator $a$, and define
$b$ similarly for $g_2$. These are primitive elliptic meridians.
Write $g_1=a^u$, $g_2=b^v$ and let $p,q$ be the orders of $a,b$.
Then $p_0\mid p$ and $q_0\mid q$. The primitive orders $p,q$
are even under either the theorem's hypothesis or the extension
stated after it. Moreover, $\langle a,b\rangle=\Gamma$, since
this subgroup contains $g_1,g_2$. The embeddings of the cyclic
factors defined by the powers $a^u,b^v$ induce an injection
$C_{p_0}*C_{q_0}\hookrightarrow C_p*C_q$ by free-product
normal forms. A nonidentity word in the kernel of the original
marked map therefore gives a nonidentity word in the kernel
of $C_p*C_q\to\Gamma$.
Theorem~\ref{thm:mainC} therefore applies to $a,b$ and gives the
models and marked representatives in Table~\ref{tab:group-models}.
The definitions in Section~\ref{sec:heckoid} identify these
models as $\mathcal L$, $\mathcal M_0$, $\mathcal M_1$, or
$\mathcal M_2$, carrying the marked weights with their edges
under the stated slope changes. Lemma~\ref{lem:labelfree} gives
$\Gamma\cong\pi_1^{\orb}(Q)$.
Proposition~\ref{thm:rigid} and
Lemma~\ref{lem:rigid-boundary-turnover} give geometric finiteness
and turnover conformal boundary. Corollary~\ref{cor:model-holonomy}
gives uniqueness of the complete holonomy.
\end{proof}

\subsection{Ends and finite volume}\label{sec:ends}

\begin{corollary}\label{cor:ends}
For a group satisfying Theorem~\ref{thm:mainC}, spherical capped
vertices are finite singular points, Euclidean capped vertices are
turnover cusps, and hyperbolic capped vertices are conformal turnover
ends. The only possible cusp sections are $S^2(2,4,4)$ and
$S^2(2,3,6)$. The group is cocompact if and only if every capped vertex
is spherical, and has finite covolume if and only if no capped vertex
is hyperbolic.

The finite-volume cases among the realized models are listed below.
The first block has even $p,q\ge4$, and the second has $p=2$ and
even $q\ge4$. Every other realized model has infinite volume.
\begin{center}\small
\begin{tabular}{@{}p{.32\textwidth}p{.23\textwidth}p{.36\textwidth}@{}}
\toprule
Model & Cocompact cases & Noncompact finite-volume cases\\
\midrule
\multicolumn{3}{@{}l@{}}{$p,q\ge4$ even}\\[4pt]
$\mathcal L(r)_{p,q}$ & Every realized model & None\\[4pt]
$\mathcal M_0(r;d)_{p,q}$ & None &
$p=q=4$, $d=2$. Two $(2,4,4)$ cusps\\[4pt]
$\mathcal M_1(r;m)_p$ &
$p=4$, $m=3$ & $p=6$, $m=3$. Two $(2,3,6)$ cusps\\[4pt]
$\mathcal M_2(r;m)_{p,q}$ & $p=q=4$, $m=3$ &
$m=3$, $p,q\in\{4,6\}$, at least one equal to $6$.
One $(2,3,6)$ cusp for each marked order $6$\\
\midrule
\multicolumn{3}{@{}l@{}}{$p=2$, $q\ge4$ even}\\[4pt]
$\mathcal L(r)_{2,q}$, two components & Every realized model & None\\[4pt]
$\mathcal M_0(r;d)_{2,q}$, knot &
$d=2$, or $(q,d)=(4,3)$ &
$(q,d)=(4,4)$. Two $(2,4,4)$ cusps.
$(q,d)=(6,3)$. Two $(2,3,6)$ cusps\\[4pt]
$\mathcal M_0(r;d)_{2,q}$, two components & None &
$q=4$, $d=2$. One $(2,4,4)$ cusp\\[4pt]
$\mathcal M_2(r;m)_{2,q}$ & $q=4$, $m=3$ &
$q=6$, $m=3$. One $(2,3,6)$ cusp\\
\bottomrule
\end{tabular}
\end{center}
\end{corollary}

\begin{proof}
Lemma~\ref{lem:degenerate} excludes infinite-weight edges and pillow
vertices. At each remaining vertex at least two incident weights are
even. Hence the Euclidean possibilities are $(2,4,4)$ and $(2,3,6)$.
Geometric finiteness makes finite covolume equivalent to empty
conformal boundary. Cocompactness also requires absence of cusps.

For $\mathcal M_0(r;d)_{p,q}$ the vertex types are $(p,q,d)$ twice
when $K(r)$ is a knot, and $(p,p,d),(q,q,d)$ otherwise.
For $\mathcal M_1(r;m)_p$ they are $(p,2,m)$ twice.
For $\mathcal M_2(r;m)_{p,q}$ they are
\[
 (p,2,2),\quad(q,2,2),\quad(p,2,m),\quad(q,2,m).
\]
The first two are spherical. With $k\ge4$ even and $m\ge3$ odd,
$(k,2,m)$ is spherical only for $(k,m)=(4,3)$ and Euclidean only
for $(k,m)=(6,3)$. With $p,q\ge4$, a vertex of type $(p,q,d)$
has reciprocal sum at most one, with equality precisely when
$(p,q,d)=(4,4,2)$. These observations give the first block.

When $p=2$, a vertex-free knot or coincident marks would force
$q=2$, so neither occurs. In the one-tunnel knot case,
$(2,q,d)$ is spherical or Euclidean precisely when
$(q-2)(d-2)\le4$, with equality in the Euclidean case.
For a one-tunnel link the vertex $(2,2,d)$ is spherical, and
$(q,q,d)$ is Euclidean only for $(q,d)=(4,2)$ and otherwise
hyperbolic. In the two-tunnel case the first three vertices
$(2,2,2),(q,2,2),(2,2,m)$ are spherical. The remaining type
$(q,2,m)$ gives the last row by the preceding calculation.
\end{proof}

\begin{example}\label{rem:triangle-example}
Let $p,q,d\ge2$ be integers with $1/p+1/q+1/d<1$.
The orientation-preserving group of the hyperbolic triangle with
angles $\pi/p,\pi/q,\pi/d$ is
\[
 \Delta(p,q,d)
 =\langle a,b\mid a^p=b^q=(ab)^d=1\rangle
 <\mathrm{PSL}_2(\mathbb R)<\PSL.
\]
Its generators $a,b$ are primitive vertex rotations of orders $p,q$.
Signed normal distance from the invariant plane identifies the quotient
with $S^2(p,q,d)\times\mathbb R$, with metric
$dt^2+\cosh^2(t)g_T$, where $g_T$ is the hyperbolic turnover metric.
In particular, it has infinite volume. Its fully capped graph is the
unknotted theta graph with weights $p,q,d$, namely
$\mathcal M_0(0;d)_{p,q}$, with the first two edges marked.
For even $p,q$ it is a model in Theorem~\ref{thm:mainC}, and its
marked order-two surgery is $\Odih(0/1;1,d)$.
\end{example}

\EvenFigureTriangle

\begin{example}[A one-cusped two-tunnel model]\label{ex:two-tunnel-lattice}
Consider the hyperbolic Coxeter tetrahedron with face reflections
$R_0,R_1,R_2,R_3$ and product orders
\[
 m_{01}=4,\qquad m_{12}=3,\qquad m_{23}=6,\qquad
 m_{02}=m_{03}=m_{13}=2.
\]
For existence, its Gram matrix has diagonal entries $1$ and
off-diagonal entries $-\cos(\pi/m_{ij})$. Its determinant is $-1/8$,
and its principal minors of order three, in the order of the omitted
indices $0,1,2,3$, are $0,1/4,1/2,1/4$. All smaller principal
minors are positive. Thus it is a finite-volume hyperbolic
tetrahedron with one ideal vertex
\cite[Sections~6.8--6.9]{Davis}.

Its orientation-preserving group $\Lambda$ is generated by the
primitive edge rotations $a=R_0R_1$ and $b=R_2R_3$ of orders $4,6$.
Indeed, put $h=R_0R_3$. Then $h^2=1$, $hah=a^{-1}$ and
$hbh=b^{-1}$, while $R_1R_2=a^{-1}hb^{-1}$.
Consequently $H=\langle a,b\rangle$ is normal in
$\Lambda=\langle a,b,h\rangle$ with index at most two.
The relation $(a^{-1}hb^{-1})^3=1$ makes the image of $h$ in
$\Lambda/H$ have order dividing three, so $H=\Lambda$.
The octahedral vertex group shows that $\Lambda\not\cong C_4*C_6$,
whose finite subgroups are cyclic.

The fully capped graph is the planar tetrahedral graph.
Remove its opposite edges of weights $2$ and $3$. The remaining
four-cycle is an unknot, with cyclic weights $4,2,6,2$, and the
removed edges lie in the two complementary discs of its containing
sphere. These are the standard upper and lower tunnels of $K(0)$.
Thus the quotient realizes $\mathcal M_2(0;3)_{4,6}$.
Its vertex types are $(2,2,4),(2,2,6),(2,3,4),(2,3,6)$, giving
exactly one cusp. Normalize the slope by $0\mapsto1$ before
using \eqref{eq:heckoid-parameter-dictionary}. Its parabolic
parent is then $S(1/2;3/2)$, so its Heckoid index is $3/2$.
\end{example}

\subsection{Consequences for groups and cusps}
\label{sec:group-cusp-consequences}

\begin{corollary}\label{cor:finite-volume-abstract}
Let $O_1,O_2$ be realized finite-volume orbifolds in
Theorem~\ref{thm:classification}, with holonomy groups
$\Gamma_1,\Gamma_2$. Every isomorphism $\Gamma_1\to\Gamma_2$ is
induced by an isometry of $\HH^3$, possibly reversing orientation.
Their abstract groups are therefore isomorphic if and only if their
unmarked fully capped labelled graphs are related by a homeomorphism
of $S^3$, with orientation reversal allowed.
\end{corollary}

\begin{proof}
For an isomorphism $\theta:\Gamma_1\to\Gamma_2$, choose a
torsion-free finite-index normal subgroup $\Gamma_1^0$ by Selberg's
lemma and put $\Gamma_2^0=\theta(\Gamma_1^0)$.
Mostow--Prasad rigidity supplies an isometry $J$ inducing
$\theta$ on $\Gamma_1^0$ \cite{Mostow,Prasad}.
For $\gamma\in\Gamma_1$, normality implies that $J\gamma J^{-1}$
and $\theta(\gamma)$ induce the same conjugation on $\Gamma_2^0$.
Their quotient is trivial, since a lattice has trivial centralizer
in the full isometry group. Indeed, a commuting isometry fixes the
dense set of endpoints of its loxodromic axes.
Thus $J$ induces $\theta$ on $\Gamma_1$.

The quotient isometry identifies singular edges and their labels.
Truncate corresponding turnover cusps and extend the boundary
maps conically over their formal caps. This gives the asserted
homeomorphism of labelled graphs in $S^3$.
Conversely, such a homeomorphism preserves the Euclidean vertex
types, so reopening them gives an orbifold homeomorphism and hence
an isomorphism of fundamental groups.
\end{proof}

A cusp translation lattice is the rank-two subgroup of translations
in its Euclidean cusp group. Its quotient gives a canonical torus
cover of the cusp section.

\begin{corollary}\label{cor:cusp-consequences}
For a group in Theorem~\ref{thm:classification}, there are at most
two cusps. If there are two, their sections have the same type.
The translation lattice of an $S^2(2,4,4)$ cusp is similar to
$\mathbb Z[i]$, and that of an $S^2(2,3,6)$ cusp is similar to
$\mathbb Z[\zeta_6]$, where $\zeta_6=e^{\pi i/3}$.
Every torus cusp section of a finite manifold cover consequently
has a basis ratio in $\mathbb Q(i)$ or $\mathbb Q(\sqrt{-3})$,
respectively.

An orbifold in Theorem~\ref{thm:classification} with an
$S^2(2,4,4)$ cusp is not commensurable with one having an
$S^2(2,3,6)$ cusp. Here commensurability means existence of a
common finite-sheeted hyperbolic orbifold cover, with orientation
reversal allowed.
\end{corollary}

\begin{proof}
Use the vertex types in the proof of Corollary~\ref{cor:ends}.
A vertex-free model has no cusp, and a one-tunnel knot has two
vertices of the same type. For a one-tunnel two-component link,
a Euclidean vertex $(k,k,d)$ with even $k\ge2$ requires
$(k,d)=(4,2)$. In a two-tunnel model the first two vertices
$(p,2,2),(q,2,2)$ are spherical. The other types $(p,2,m),(q,2,m)$,
where $m\ge3$ is odd, are Euclidean only when the marked order is
$6$ and $m=3$. This proves the cusp count and uniformity,
including in infinite volume and when $p=2$.

For the square or hexagonal cusp, its translation lattice $L$
is invariant under multiplication by $\epsilon=i$ or $\zeta_6$,
respectively. Choose $0\ne\alpha\in L$. Then
$\alpha\mathbb Z[\epsilon]$ has finite index in $L$, so
$\alpha^{-1}L$ is a fractional ideal of $\mathbb Z[\epsilon]$.
Both rings are Euclidean. Their fractional ideals are principal,
which gives the stated similarity classes.
Finite manifold covers give finite-index translation sublattices,
whose basis ratios lie in the corresponding quadratic field.
Such sublattices need not themselves be square or hexagonal
\cite[Section~2.2 and discussion preceding Theorem~3.7]{LongReid}.

Suppose the two orbifolds had a common finite-sheeted hyperbolic
orbifold cover. Pass to a further manifold cover by Selberg's lemma
and choose a torus cusp above the square cusp. Its rank-two
parabolic stabilizer has finite index in its stabilizers in both
base groups, so it also covers a cusp on the other side.
Uniformity makes that cusp hexagonal. A basis ratio of the torus
lattice would then lie in
$\mathbb Q(i)\cap\mathbb Q(\sqrt{-3})=\mathbb Q$, contradicting
its nonzero imaginary part. Complex conjugation preserves both
fields, so orientation reversal gives the same contradiction.
\end{proof}

\begin{corollary}\label{cor:no-knot-cover}
No orbifold in Theorem~\ref{thm:classification} is finitely covered
by a hyperbolic knot complement in $S^3$.
\end{corollary}

\begin{proof}
If a hyperbolic knot complement finitely covers $O$, then $O$ has
finite volume and exactly one cusp. Its section is $S^2(2,4,4)$
or $S^2(2,3,6)$ by Corollary~\ref{cor:ends}.
The first is excluded by \cite[Theorem~1.1]{HoffmanCusp}.
In the second case, \cite[Theorem~1.2]{HoffmanCusp} gives
$\Gamma^{\mathrm{ab}}\cong\mathbb Z/2\mathbb Z$.

The end table forces a two-tunnel model of weights $2,3$.
Its even subgraph is $K(r)\cup\tau_+$, with four vertices and
five edges. It is connected because $K(r)$ is connected in the
knot case and the upper tunnel joins the two components otherwise.
Lemma~\ref{lem:mod2weighted} therefore gives
\[
 \dim_{\mathbb F_2}H_1(\Gamma;\mathbb F_2)=5-4+1=2,
\]
contradicting $\Gamma^{\mathrm{ab}}\cong\mathbb Z/2\mathbb Z$.
\end{proof}

\section{The limitation of the parity hypothesis}\label{sec:parity-limits}

The divisibilities $2\mid p$ and $2\mid q$ make the order-two
weight change induce the quotient
$\Gamma\twoheadrightarrow\Gamma/\nc{a^2,b^2}$
in Lemma~\ref{lem:paritysurgery}. If a marked order is odd,
adding the corresponding square relation kills that generator.
The following examples show that the spatial conclusion can
also fail when an odd marked order is allowed.

\begin{proposition}\label{prop:odd-counterexamples}
For $n=4,5$ there is a cocompact Kleinian group $\Gamma_n$, generated
by primitive elliptics of orders $3,n$, which is not the free product
$C_3*C_n$. Its quotient has underlying space $S^3$ and singular graph
the planar tetrahedral graph. The four vertex groups are
$A_4,A_4,S_4,S_4$ for $n=4$, and $A_4,A_4,A_5,A_5$ for $n=5$.
In particular these quotients are not the canonical spatial Heckoid
models of Section~\ref{sec:heckoid}, even if arbitrary finite marked
orders are allowed there.
\end{proposition}

\begin{proof}
Put $t=2\cos(\pi/n)$ and consider the Gram matrix
\[
 B=\begin{pmatrix}
 1&-1/2&0&-t/2\\
 -1/2&1&-1/2&0\\
 0&-1/2&1&-1/2\\
 -t/2&0&-1/2&1
 \end{pmatrix}.
\]
Its determinant is $-(t-1)(3t+5)/16<0$. Its principal minors of
order three are $1/2,(3-t^2)/4,(3-t^2)/4,1/2$, all positive, and
the minors of orders one and two are positive. Thus $B$ has signature
$(3,1)$ and defines a compact hyperbolic Coxeter tetrahedron
\cite[Theorem~6.9.1 and Section~6.12]{Davis}. Let $R_0,\ldots,R_3$
be its face reflections. Their product orders are
\[
 m_{01}=m_{12}=m_{23}=3,\qquad m_{03}=n,\qquad
 m_{02}=m_{13}=2.
\]
The orientation-preserving subgroup $\Gamma_n$ is generated by
$x=R_0R_1$, $f=R_1R_2$ and $y=R_2R_3$. Put $g=R_0R_3=xfy$.
The following identities hold:
\begin{align*}
 x&=W_4(f,g)=g^2fgf^{-1}g^{-1}fgfg^{-1}f^{-1}g^2
 &&(b=4),\\
 x&=W_5(f,g)=fgfg^{-1}f^{-1}gfgf^{-1}g^{-1}fgfg^{-1}f^{-1}gf
 &&(b=5).
\end{align*}
They can be checked in the faithful reflection representation
$R_i(v)=v-2B(v,e_i)e_i$, using respectively $t^2=2$ and $t^2=t+1$.
Since $y=f^{-1}W_b(f,g)^{-1}g$, the pair $f,g$ generates $\Gamma_n$.
Products of adjacent face reflections are once-around edge rotations,
so these generators are primitive of orders $3,n$.

The quotient is the double of the tetrahedron along its faces. Its
singular graph is the tetrahedral one-skeleton, with vertex types
$(2,3,3)$ twice and $(2,3,n)$ twice. The associated rotation groups
are exactly those in the statement. Their noncyclic finite subgroups
exclude a free-product isomorphism, since every finite subgroup of
$C_3*C_n$ is conjugate into a cyclic factor \cite{Serre}.

A link model has no vertices and a one-tunnel model has at most two.
A two-tunnel canonical model has two actual dihedral vertices
of types $(2,2,p)$ and $(2,2,q)$, at the feet of its weight-two
tunnel. Our four vertices are all actual and non-dihedral. Hence none of these weighted spatial models is the
quotient just constructed.
\end{proof}

For $n=4$, the lattice $\Gamma_4$ is the orientation-preserving
subgroup of the tetrahedral reflection group
$T_5[2,3,3;2,3,4]$, identified in
\cite[Section~1.3, no.~8 and the paragraph following the table]{MSY}.
Its four non-dihedral vertices give an exception to the list of
singular graphs accompanying \cite[Figure~1]{MSY}.

\end{document}